\documentclass{amsart}
\usepackage{amsthm,amsfonts,amsmath,amscd,amssymb,latexsym,epsfig,eucal,float}
\usepackage{multirow}

\newtheorem{theorem}{Theorem}

\def\BOne{{\mathchoice {\rm 1\mskip-4mu l} {\rm 1\mskip-4mu l}
                          {\rm 1\mskip-4.5mu l} {\rm 1\mskip-5mu l}}}

\def\<{\langle}

\def\>{\rangle}

\def\mb#1{{\mathbb #1}}
\def\mc#1{{\mathcal #1}}

\begin{document}

\title[The vibrations of thin plates]{The vibrations of thin plates}
\author{Santiago R. Simanca}
\thanks{Supported by the Simons Foundation Visiting Professorship award 
number 657746.}
\address{Department of Mathematics, Courant Institute of Mathematical Sciences,
251 Mercer St., New York, NY 10012} 
\email{srs2@cims.nyu.edu}

\begin{abstract}
We describe the equations of motion of an incompressible elastic 
body $\Omega$ in 3-space acted on by an external pressure force, and the Newton
iteration scheme that proves the well-posedness of the resulting 
initial value problem for its equations of motion on $C^{k,\alpha}$ spaces. 
We use the first iterate 
of this Newton scheme as an approximation to the actual vibration motion of 
the body, and given a (finite) triangulation $K$ of it, produce an 
algorithm that computes it, employing the direct sum of the space of PL 
vector fields associated to the oriented edges and faces of the 
first barycentric subdivision $K'$ of $K$ (the metric duals of the 
Whitney forms of $K'$ in degree one, and the metric duals of the 
local Hodge $*$ of the Whitney forms in degree two, respectively) as 
the discretizing space. 
These vector fields, which capture the algebraic topology properties of 
$\Omega$, encode them into the solution of the weak version of the 
linearized equations of motion about a stationary point, the essential component
in the finding of the first iterate in the alluded Newton scheme. This
allows for the selection of appropriate choices of $K$,
relative to the geometry of $\Omega$, for which the algorithm produces 
solutions that accurately describe the vibration of 
thin plates in a computationally efficient manner. We use
these to study the resonance modes of the vibration of these plates, and 
carry out several relevant simulations, the results of which are 
all consistent with known vibration patterns of thin plates derived 
experimentally.
\end{abstract}

\maketitle 

\section{Introduction} \label{s1}
The motion of an incompressible elastodynamic body $\Omega$ 
is described by a path of embeddings 
$t\rightarrow \eta(t): \Omega \hookrightarrow \mb{R}^3$ 
that satisfies a nonlinear pseudodifferential wave equation, and
with the boundary $\partial (\eta(t)\Omega)$, which is free to move, doing so
following some conditions in the normal directions.  
The spatial component of the wave equation is an elliptic operator 
determined by a tensor $W$, which encodes the internal energy stored in 
$\Omega$ at the microscopic level as it is deformed in the various directions,
exactly as a linear spring stores energy when it is compressed, or elongated.
This elliptic operator has a nonlocal part, a correction term introduced by
the gradient of a pressure function, that ensures that the motion stays 
incompressible at all time (that is to say, volume preserving at the
infinitesimal level everywhere). And since $\eta(t)$ maps points on 
$\partial \Omega$ to points on $\partial \eta(t) \Omega$, any tangential 
change over the boundary must be compensated for by a corresponding 
change in the normal direction, so that the incompressible condition holds at 
those points as well. The mechanism by which this boundary motion happens
is thus, a function of the stored energy tensor $W$ also.

At least for a short time, the initial value Cauchy problem for this
nonlinear pseudodifferential wave is well-posed \cite{ebsi2}.
All particle points that are 
deformed a sufficiently small amount tend to go back to their equilibrium
state, much like the spring does while it is deformed in its elastic regime.
The pseudodifferential terms in the equation make the entire body feel 
these deformations at one point instantly anywhere else in the body, but they 
are initially so tiny that their effect on the nonlinear terms of the 
equation are negligible, and the body moves then as if its motion were
being ruled by a differential linear wave equation instead.
While the body keeps moving, eventually, the effect of these tiny local
changes may add up to a point where the nonlinear terms in the equation, local
and nonlocal, could enhance the effect they produce on the overall motion, 
making them no longer negligible as they were at the beginning.  
The body could then become irreparably deformed at locations
where the motion gets driven into the plastic regime, developing cracks by 
inelastic shearing, or breakages by inelastic pull, 
if the effects of the tiny deformations grow to be 
sufficiently large that the nonlinear terms in the equation become the most
significant, and quite large, at those locations where the crack or break 
is occurring. Up until the moment when singularities develop, if at all, 
the boundary moves so that the {\it directional derivative} of
$W$ along the exterior normal $N$ of $\partial \Omega$, at a boundary point
$x$, is a vector field that points in the direction of
the exterior normal $\nu$ of $\partial(\eta(t) \Omega)$, at $\eta(t)(x)$. 

The condition ruling the motion of the boundary makes visible the significant 
additional challenge in the study of the motion of very thin $\Omega$s,  
bounded three dimensional bodies with one of the dimensions at least
one order of magnitude smaller in length than the other two, thus, 
geometrically, $3$d bodies that almost degenerate into $2$d plates.
We have far apart pairs of boundary points on ``oppossite sides'' of a
thin plate that are separated by a very small distance within the plate. 
At each of the
points in these pairs, the exterior normals to the boundary point 
in directions almost opposite to each other, and so, while the motion does not
develop singularities, these boundary points are being pulled further apart,
elongating locally the body in the thin direction, or compressed into each 
other, further thinning the body at location. This phenomena 
accelerates the plausible formation of singularities in the motion of the 
body, a direct consequence of its quasi geometric degeneration.

An stationary thin plate is caused to vibrate when acted on by an external 
periodic pressure force; equivalently, a thin plate that moves uniformly 
through space
is caused to vibrate by the action of the air pressure on it when the air 
pressure in the area where the plate is moving changes (somewhat) periodically.
At certain frequencies of the external pressure force, the 
body responds and vibrates by resonance, the nodal and antinodal 
configuration points of these waves characteristic of the plate
at its eigenfrequency resonance modes. Several acoustic experiments 
serve to illustrate this situation, notably, those carried out by
F\'elix Savart as far back as 1830, and which were built on a method 
developed by Ernst Chladni (see, \cite[$2^{\rm nd}$ 
column, page 171]{hu}; this reference describes various other types of 
related acoustic experiments also).  
Today, particularly among luthiers, the nodal lines of plates vibrating 
at these frequencies are called the Chladni lines, and the entire 
portrait that they make is called the Chladni pattern.    

In this article, we study the theoretical, and practical foundations of these
situations. We assume that the motion of $\Omega$ is incompressible, and 
proceeding in general, describe firstly the equations that rule the motion in 
analogous circumstances. We then produce an algorithm to compute a judicious 
numerical solution approximation to the motion, and carry out various 
simulations with it, displaying the resulting 
Chladni patterns for a handful of vibrating thin $\Omega$s.

In our simulations, we suppose that the plates are  
made of {\it orthotropic} elastic Hooke material, 
wood to be specific. The tensor $W$ of these bodies 
is characterized by nine independent parameters, which,
together with the density, determine the equation of motion. For thin 
$\Omega$s of four distinct relevant geometries, assuming further that the
parameters are constant throughout the body, we compute numerical 
solution approximations of the equations, and display the Chladni patterns for
each of these plates vibrating by resonance at five 
different frequencies, these frequencies chosen to be those considered as 
the most important 
vibration modes in the tuning of a violin plate.  
We do so with an eye towards comparing, and validating our results, against
those derived in the acoustic experiments mentioned above. 

The finding of accurate numerical solutions to the equations of motion of 
incompressible elastic bodies is a matter of interest in its own right, 
and quite a difficult problem in general given the pseudodifferential nature 
of the nonlinear equation being solved, with changes to the 
solution in a neighborhood of any point affecting its value everywhere else 
all at once. These difficulties are further enlarged if
we deal with bodies that are almost degenerate, very thin in one direction.
We overcome both of these difficulties by implementing two key ideas 
that arise after taking a close look at the method of proof of the 
well-posedness of the equation of motion.

The said proof is based on the contraction mapping principle, 
originally carried out working on Sobolev spaces \cite{ebsi2}, and later
on extended and shown to work on $C^{k,\alpha}$ spaces as well \cite{sim}.
The benefit of the latter method is twofold: On the one 
hand, we get optimal regularity results for the analysis of the
Cauchy problem, which can be started if we assume merely that $\eta(t)$ is
a $C^{2,\alpha}$ curve of embeddings; and on the other hand, the
$C^{k,\alpha}$-spaces are better suited to analyze the question of 
consistency of any numerical solution of the equations of motion 
that we might propose.
Our algorithm computes numerically the weak solution of the linearized
equation that yields the first orbit point in the Newton iteration scheme 
used to prove the theorem. Since this orbit point lies 
in $C^{1,\alpha}$, the weak solution of the equation can be construed as an 
element of $L^2$ that has weak derivatives in $L^2$ also. 
By the topological nature of the unknown in the equation, it is natural
to resolve this problem introducing a finite triangulation $K$ in the body. 
We may then use, as discretizing spaces for this weak solution, the 
direct sum of the spaces of PL vector
fields given by the metric duals of the degree one Whitney forms, and the 
metric duals of the Hodge $*$ of the degree two Whitney forms of the 
first barycentric subdivision $K'$ of $K$, respectively. 
A good such choice of $K$ allows us to overcome, as efficiently as possible,  
the two difficulties above inherent in the problem.

Indeed, as this weak solution yields an approximation to the actual solution of
the linearized pseudodifferential equations of motion, whose velocity field 
is a divergence-free vector field, the discretizing spaces for its numerical 
casting should have encoded into them the algebraic topology of gradient 
and divergence-free fields.  The weak solution is not, in itself, 
divergence-free, and has a gradient component also, the latter being 
relatively small given the compatible 
initial conditions for the Cauchy problem that it satisfies. 
The metric duals of the Whitney forms of $K'$ in degrees less or equal than
one, or their Hodge $*$ in degrees greater or equal than two,
produce a simplicial complex with functions as the Abelian groups of 
the complex in degrees zero, and three, and vector fields in degrees 
one, and two, respectively. (Notice the important fact here that $K'$ is a 
naturally oriented simplicial complex.)  
The cohomology of this complex is the cohomology of the body, 
and in degrees one and two, the cycles are the Abelian subgroups
of PL gradient fields associated to the edges, and 
divergence-free vector fields associated to the faces, respectively. The 
sought after solution is discretized as an element of the direct sum 
of the cochain groups of this complex in these two degrees. The closer this
solution gets to be the sum of an actual gradient and actual 
divergence-free field, 
the closer it will get to be an element of the subspace given by the 
direct sum of the alluded cycle subspaces. The vector of coefficients of 
the linear combination producing the discretized solution is found by 
solving the linear system of second 
order ordinary differential equations arising from the weak formulation of 
the equation over this discretizing space, a square system of size equal to 
the sum of the number of edges and faces in $K'$.  
The global nature of the pseudodifferential wave equation is thus transformed
into a local problem for a rather large, but very tractable, linear system
of differential equations.   

The resonance vibration patterns that we intent to describe involve 
primarily small vibrations of the bodies under consideration, and these
can be approximated well by the solution to the linearized equations of 
motion about the canonical stationary state, which in turn is described by
our numerical solution wave. With this in mind, we see that the problem 
generated by the effects of the almost degenerate nature of a thin plate on 
our numerical algorithm is resolved by choosing $K$ to have sufficiently many 
simplices, each of them with an aspect ratio in the order of one, and so, 
a fortiori, a triangulation $K$ with a large number of simplices in it. 
For then we have the necessary resolution for the numerical approximation 
to the motion to be accurate at any point of the plate. 
The uniformly distributed {\it oriented} edges and faces of $K'$ that ensues 
results into an almost constant number of these being placed at every location 
of the plate. At the local level, the Whitney vector fields associated to a
face and its bounding edges interact with one another, and the interaction 
spills over to neighbor faces and edges, overall producing a global 
interaction of all of 
the Whitney PL vector fields with each other. This allows for
the numerical solution to feel at any point the contributions to the 
vibration modes arising from all the points in the 
body, including the many far apart boundary points that are separated by a 
very small distance within the body, with the accuracy of the approximation
improving as we enlarge the local almost constant number of edges and faces.  
We pay a larger computational price the larger we choose this local constant 
number to be, triangulating the body with the appropriate resolution,
but the accuracy of the results increases as we do so. We are
able to determine an appropriate resolution here (relative to 
the thinness of the plate) leading to a satisfying accuracy, and 
manage the computational complexity of the problem
with this choice of resolution for $K$ using very modest resources.
   
All the known type of waves within the body, Lamb, Rayleigh, shear, or 
otherwise, fall within a single framework. They result from the 
elastic interaction of the material points that compose it, whose potential
energy is codified into the tensor $W$.
The solutions we construct numerically to describe these waves are built 
to maintain the global topological constraint imposed by the 
incompressibility condition in time, yielding accurate approximations to 
the actual motion of the body. The innovative use of our discretizing 
spaces counterbalances the need to go through the computational complexities 
inherent in the problem, and in exchange we are able to produce results 
that are faithful to the physical reality of the motion while in the 
elastic regime.

In our simulations for orthotropic bodies, the tensor $W$ is assumed to
be covariantly constant. However, our approach works as well in the study 
of the vibrations of incompressible bodies whose stored energy tensor $W$ is 
(or is assumed to be) just $C^2$ differentiable, as well as for bodies 
that have portions of the boundary fixed, while the rest remains free to move, 
or incompressible bodies of this more general type (even perfect fluids for 
that matter) embedded into nonorientable $3$d Riemannian spaces, instead of 
Euclidean $\mb{R}^3$. A case of particular interest would be the treatment 
of the problem for isotropic functionally graded plates \cite{mahi}; 
the smaller number of elastic parameters would make that treatment far easier 
by comparison, even if the values of the parameters now vary across the body. 

Our simulations are computationally more complex than 
approaches describing aspects of the motion in terms of an ad hoc small number 
of degrees of freedom, and basic assumptions, but the physical meaning of 
their results cannot be questioned.

All of our simulations are in close correspondence with those carried out in 
\cite{sim}, and fit well with the alluded classical experiments on violin 
plates above, as did the ones before. But we now expand significantly on 
accuracy, and range of applicability, the 
reasons why to be clarified in detail below, when we get the opportunity to 
contrast the equations solved, and the manner in which they were solved then, 
and are solved now. We use a simplifying analogy here that, loosely 
speaking, conveys 
quickly to the reader the differences in these works, and why the results 
are good 
in the regime where they are both applicable. For suppose that we look at 
a nonlinear, and nonhomogeneous system of ordinary differential equations, and 
approximate its solution with trivial data using the first iterate of 
Euler's method, or approximate its solution with some arbitrary data using 
the first iterate of the Cauchy-Peano method. The numerical approximation to 
the solution of the equations of motion of incompressible bodies in 
\cite{sim} is to the former of these approaches, what the numerical 
solution of this same equation here is to the latter.

\subsection{Organization of the article}
In \S\ref{s2} we state the equations of motion of incompressible 
elastodynamic bodies, and briefly sketch their slight modifications leading to
the proof of the well-posedness of the Cauchy problem \cite{ebsi2,sim}. We 
emphasize the nonhomogeneous version of the equation, restate its 
linearization at an 
arbitrary point, and particularize the latter at a starting point 
of the form $(\eta(0),\dot{\eta}(0))=(\BOne, w)$, $w$ a divergence-free
given field. We then give a full description of the first iterate in the
Newton scheme that proves the well-posedness (in terms of $w$, and the
nonhomogeneous term in the equation). In \S\ref{s3}, we recall 
the notion of a generalized Hooke body, its stress and strain tensors, 
and the particular properties of an orthotropic one, 
together with the nine elastic parameters that characterize its 
stored energy tensor $W$. We show also the values of these parameters
for the orthotropic material that are used in our simulations later on, in 
\S\ref{s6}. In \S\ref{s4}, we describe the Whitney forms of a $3$d manifold 
with boundary $M$, triangulated by the barycentric subdivision $K'$ of a 
finite triangulation $K$, and the simplicial complex they give rise to, whose 
homology is the cohomology of $M$. The direct sum of the groups of this
complex in degrees one and two, the Whitney vector fields associated to the
edges and faces of $K'$, is shown to parallel the usual decomposition of
a vector field into a gradient field, and a divergence-free field, which are
$L^2$ orthogonal to each other, a fact at the heart of the Poincar\'e duality 
for the simplicial complex, and which makes of this direct sum space the 
natural choice to discretize the weak solution that we pursue numerically,
given its algebraic {\it and} geometric content.  The algorithm
to solve numerically the weak solution of the linearization of the 
modified equation of motion used to prove the well-posedness is
explained in detail in \S\ref{s5}, as well as an algorithm that 
from it, computes the full fledged first iterate of the Newton scheme 
employed in the proof. 
In \S\ref{s6}, we present the simulations, 
and some additional numerical computations to support the virtues of our 
approach, in each case discussing the numerical complexity. We contrast 
the results obtained against known experiments, as a way of validating them.
We end with some remarks of interest in \S\ref{s7}, synthetizing the essence
of the proofs of various results in the article, and pointing 
towards generalizations of various aspects of our work here.

\section{Incompressible motions}
\label{s2}
We begin by summarizing the basics of the equations of motion 
of elastodynamic incompressible bounded bodies, and the essence of the argument
that treats the well-posedness of the associated free-boundary initial value 
problem. We work with three dimensional bodies embedded in $\mb{R}^3$,
though the results extend to any dimension $n$.
We then bridge this to the nonhomogeneous problem resulting from
the motion of the body under the influence of an external force, and
in that context, we describe the fixed point iteration scheme argument that
leads to the well-posedness of the equations of motion for a short time. We
give special attention to the first complete iterate of this Newton scheme 
when the body starts its $F$ driven motion from the rest position, as we
shall use it as an approximation to the actual motion of 
the body, which is given by the fixed point of the scheme instead.
The reader may want to consult \cite{ebsi,ebsi2}, and \cite{sim}, for relevant 
details on both topics.

We let $\Omega$ be a bounded domain in $\mb{R}^3$, whose  
boundary is of class $C^{1,\alpha}$. We assume
that the density $\rho$ of the material filling $\Omega$ is 
constant. The motion of $\Omega$ is encoded into a curve
$$
\eta(t) : \Omega \hookrightarrow \mb{R}^3
$$
of embeddings of $\Omega $ into $\mb{R}^3$ that preserve 
volume. The path $t \rightarrow \eta (t)(x)$ denotes the position at 
time $t$ of a particle initially at $x\in \Omega$. We denote by
$D \eta(t)(x)= (\partial_{x^i}\eta^j(t)(x))$ the deformation gradient.

The material properties of $\Omega$ are characterized by its
stored energy function $W$, which it is assumed
to be a function of the deformation gradient, $W=W(D\eta)$. (This function 
is the quadratic form associated associated to the tensor $W$, the 
reason why we shall refer to both, the function and the tensor, 
using the same term, see \S \ref{s3} below). In  
the presence of no external forces, the trajectory of the body
is an extremal path of the Lagrangian
\begin{equation} \label{lag}
{\mc{L}}(\eta )=\frac{1}{2}\int _{0}^{T}\int _{\Omega }
\rho \| \dot{\eta } (t)(x)\| ^{2}dxdt - \int _{0}^{T}\int _{\Omega }
W(\partial _{i} \eta ^{\alpha }(t)(x))dxdt\; ,
\end{equation}
where the stationary points of $\mc{L}(\eta)$ are searched for among 
incompressible variations of $\eta(t)$. The motion 
is described by the solution to the system of equations
\begin{equation} \label{eq1}
J(\eta(t))=\det{D \eta (t)(x)}= 1\, ,
\end{equation}
\begin{equation}\label{eq2}
\rho \ddot{\eta }(t)(x)-{\rm Div}\, W^{'}(D\eta (t)(x))=
\nabla p(t)(\eta (t)(x)) \, ,
\end{equation}
\begin{equation}
W^{'}(D\eta )N+p(t)(\eta(t)(x))J^{b}(\eta)\nu \circ \eta =0 \; {\rm on}
\; \partial \Omega \, ,
\label{eq3}
\end{equation}
where the pressure function $p(t): \eta (t)(\Omega ) \rightarrow \mb{R}$ 
is a pseudodifferential operator in $(\eta, \dot{\eta})$.   
Here, $W^{'}$ is the derivative of $W$ with respect to the 
variables $D\eta $, ${\rm Div}$ is the divergence of $W^{'}$ with respect to
the material coordinates $x$, 
$N$ and $\nu$ are the
unit vectors normal to $\partial \Omega$ and $\partial (\eta(t)(\Omega))$,
respectively, and $J^b(\eta)$ is the Jacobian determinant of $\eta$
restricted to the boundary.

In coordinates $\eta^\alpha=\eta ^{\alpha }(x^{1}, x^2 , x^{3})$, we have that
$$ 
({\rm Div}\,  W^{'})^{\alpha }
=\partial _{x^{i}}\left( \frac{\partial W}{\partial (\partial_{i}\eta ^{\alpha
}
)}\right)\stackrel{{\rm def}}{=} A_{ij}^{\alpha \beta }(D\eta )\partial _{i}
\partial _{j}\eta ^{\beta }
$$ 
and
$$ 
(W^{'}(D\eta )N)^{\alpha }=
\frac{\partial W}{\partial (
\partial_{i}\eta ^{\alpha })}N^{i}\, , 
$$
respectively, and the system above is given by 
$$
\rho \ddot{\eta }^{\alpha }  =  A_{ij}^{\alpha \beta }(D\eta )\partial _{i}
\partial _{j}\eta ^{\beta }+(\partial _{\alpha }(p(t)))\circ \eta \, , 
$$
$$
\frac{\partial W}{\partial (
\partial_{i}\eta ^{\alpha })}N^{i} +p(t)(\eta(t)(x))J^b(\eta) (\nu\circ \eta)^
\alpha  = 0 \, . 
$$
We require the stored energy function $W$ to be {\it coercive}, so we assume 
that the operator $A_{ij}^{\alpha \beta}(D\eta)\partial_i \partial_j$ is
uniformly elliptic in a neighborhood of the curve $\eta(t)(x)$.

Given an operator $F$, we define the operator $F_{\eta}$ by 
$F_{\eta}u=(F(u\circ \eta^{-1})) \circ \eta$. If we now differentiate 
(\ref{eq1}) with respect to $t$, we obtain that 
${\rm div}_{\eta}\dot{\eta}=0$, and upon a second differentiation, we have that
$$
{\rm div}_{\eta} \ddot{\eta}= -[ (\dot{\eta}\circ \eta^{-1})
\cdot \nabla, {\rm div}]_{\eta}
\dot{\eta}={\rm trace} (D_{\eta}\dot{\eta})^2 \, .
$$
Since embeddings map boundary points to boundary points, 
we have that $W'(D\eta)N$ is perpendicular to $\partial_T \eta$ 
(where $T$ is any vector tangent to $\partial \Omega$), and the motion is  
described by the equivalent first order system 
\begin{equation}
\frac{d}{dt}\left( \begin{array}{c}
                   \eta \\ \rho \dot{\eta } 
                   \end{array}\right) = \left(
                                            \begin{array}{c}
                                            \dot{\eta } \\ A(\eta ,\dot{\eta })
                                        \end{array}\right) 
        \stackrel{{\rm def}}{=} F(\eta ,\dot{\eta }) \; ,
\label{eqs} 
\end{equation}
where
$$
A(\eta ,\dot{\eta }) = {\rm Div}\, W'(D\eta) + \nabla _{\eta }q \, ,
$$ 
and where the pressure function $q=p\circ \eta$ solves the boundary value 
problem
\begin{equation}
\begin{array}{c@{=}l}
\Delta _{\eta }q \; & \; -{\rm div}_{\eta }{\rm Div}\, W'(D\eta) +\rho \, 
{\rm trace}(D_{\eta } \dot{\eta })^{2} \vspace{1mm} \\
q\mid _{\partial \Omega }\; & \;-{\displaystyle \frac{\langle W'(D\eta) N
 , \nu \circ \eta  \rangle }{J^{b}(\eta )}} \; .
\end{array}
\label{eqb} 
\end{equation}

By the hypothesis on $W$, the boundary value problem (\ref{eqb}) is elliptic, 
and has a unique pressure function solution $q=q(\eta,\dot{\eta})$, which
is a nonlocal pseudodifferential operator in $(\eta, \dot{\eta})$.
If $w$ is a constant divergence-free field, the system (\ref{eqs}) admits 
the time independent curve $(\eta,\dot{\eta})=(\BOne,w)$ as a solution, 
with pressure function the constant $q=-\< W'(\BOne)N,N\>$. 

\begin{theorem}
{\rm (\cite[Theorem 5.53 and \S 6]{ebsi2}, \cite[Theorem 2]{sim})}
Under the hypothesis above on $\Omega$ and $W$,
the Cauchy problem for {\rm (\ref{eqs})} with 
initial condition $(\eta (0),\dot{\eta}(0))=(\BOne, w)$, 
$w$ a divergence-free vector field, is well-posed over a time interval 
whose length depends only upon a suitable norm of the Cauchy data.
\end{theorem}

The two proofs of this result use a contraction mapping principle working
on Sobolev spaces of sufficiently high order \cite{ebsi2}, or in
 $C^{k,\alpha}$ spaces with $k\geq 2$ \cite{sim}, respectively. The  
technical difficulties imposed by condition (\ref{eq1}) are overcome 
 by modifying equations (\ref{eq2}) and 
(\ref{eq3}) slightly, and considering instead the equation and 
boundary conditions   
\begin{equation} \label{new2}
\begin{array}{c}
\rho \ddot{\eta }(t)(x)=A_{\lambda}(\eta, \dot{\eta}) 
 \stackrel{def}{=}{\rm Div}\, W^{'}(D\eta (t)(x))+
J(\eta(t))\nabla_\eta q +\lambda J(\eta(t)) \nabla_\eta 
J(\eta(t) ) \, , \vspace{1mm}\\
\text{$\< W'(D\eta ) N , \partial_T \eta\>=0$, \quad $J(\eta)=1$ \quad 
on $\partial \Omega $,}
\end{array} 
\end{equation}
$\lambda$ some positive constant chosen, and fixed a priori.
Here $T$ is any vector tangent to $\partial \Omega$, and the scalar 
function $q$ solves the boundary value problem
\begin{equation} \label{newq}
\begin{array}{rcl}
L(\eta)q \stackrel{def}{=}{\rm div}_{\eta} J(\eta) \nabla_{\eta }q 
& = & -{\rm div}_{\eta }{\rm Div}\, W'(D\eta) +\rho \,
{\rm trace}(D_{\eta } \dot{\eta })^{2}\, ,  \vspace{1mm} \\
q\mid _{\partial \Omega }\; & = &  \;-{\displaystyle \frac{\langle W'(D\eta) N
 , \nu \circ \eta  \rangle }{J^{b}(\eta )}} \; . 
\end{array}
\end{equation}
We write this equation as the first order system 
\begin{equation} \label{mod}
\frac{d}{dt}\left( \begin{array}{c}
                   \eta \\
                   \rho \dot{\eta}
\end{array} \right)
= \left( \begin{array}{c}
                   \dot{\eta} \\
                   A_{\lambda}(\eta,\dot{\eta})
\end{array} \right) \stackrel{def}{=}G_\lambda(\eta,\dot{\eta}) \, .
\end{equation}
A solution $(\eta(t),\dot{\eta}(t))$ of this system for which
$\eta(t)$ satisfies (\ref{eq1}), is a solution of 
(\ref{eqs}). And vice versa.

The linearization of (\ref{mod}) at $(\eta, \dot{\eta})$,  
in the direction of $(u,v)$ yields the system     
\begin{equation} \label{lal}
\frac{d}{dt}\left( \begin{array}{c}
                   u \\ \rho v 
                   \end{array}\right)  
=D_{(\eta , \dot{\eta })}G_{\lambda} \left( \begin{array}{c}
                                 u \\ v 
                                 \end{array}
                          \right) = 
      \left(
           \begin{array}{c}
           v \\ D_{(\eta , \dot{\eta })}A_{\lambda} 
                  \left( 
                       \begin{array}{c}
                       u \\ v 
                       \end{array}
                  \right) 
           \end{array}
      \right) \, .   
\end{equation}
We let $(u(t),v(t))$ be its solution with Cauchy data compatible with 
the given Cauchy data $(\BOne, w)$ for (\ref{eqs}).
Then, for sufficiently large $\gamma$, we solve the equation
\begin{equation} \label{stl}
\gamma \left( \begin{array}{c}
                   \zeta \\
                   \rho \dot{\zeta}
\end{array} \right)
- G_\lambda(\zeta,\dot{\zeta}) =
-\left( \begin{array}{c}
                   u(t) \\
                   \rho v(t) 
\end{array} \right) +\gamma\left( \left( \begin{array}{c}
                   \BOne \\
                   \rho w 
\end{array} \right)
+\int_0^t \left( \begin{array}{c}
                   u(s) \\
                   \rho v(s) 
\end{array} \right) ds\right)  \, ,
\end{equation}
for $(\zeta,\dot{\zeta})$, for each fixed $t$. There results a mapping
$$
\mc{M}: (\eta(t), \dot{\eta}(t)) \rightarrow (\zeta(t),\dot{\zeta}(t))
$$
that, over a suitable domain of curves defined on some time interval, 
is a contraction. Its fixed point $(\eta,\dot{\eta})$ solves (\ref{mod}), 
and the diffeomorphism $\eta(t)$ so produced is volume preserving, and 
(\ref{eq1}) holds. Thus, $(\eta(t),\dot{\eta}(t))$ is the desired solution 
of (\ref{eqs}) with the said Cauchy data, 
 and over the time interval
where it is defined, $\eta(t)$ depends continuously upon the initial conditions.

We assume now that the body $\Omega$ is acted on by an  external pressure 
force $F$, and rederive the  
nonhomogeneous version of the approach above to well-posedness.
Thus, starting with pairs $(\eta, \dot{\eta})$ such that 
$(\eta(0),\dot{\eta} (0))=(\BOne, w)$, we solve the nonhomogeneous version 
of (\ref{lal}) given by
\begin{equation} \label{lah}
\frac{d}{dt}\left( \begin{array}{c}
                   u \\ \rho v
                   \end{array}\right)
=D_{(\eta , \dot{\eta })}G_{\lambda} \left( \begin{array}{c}
                                 u \\ v
                                 \end{array}
                          \right) =
      \left(
           \begin{array}{c}
           v \\ D_{(\eta , \dot{\eta })}A_{\lambda}
                  \left(
                       \begin{array}{c}
                       u \\ v
                       \end{array}
                  \right)
           \end{array}
      \right) +\left( \begin{array}{c} 
 0 \\ F
\end{array} \right) 
\end{equation}
with compatible Cauchy data. If $(u,v)$ is the solution, we then 
consider the equation
\begin{equation} \label{imp}
\gamma \left( \begin{array}{c}
                   \zeta \\
                   \rho \dot{\zeta}
\end{array} \right)
- G_\lambda(\zeta,\dot{\zeta}) =
-\left( \begin{array}{c}
                   u(t) \\
                   \rho v(t)
\end{array} \right) +
\int_0^t  \left( \begin{array}{c}
                   0  \\
                   F(s)
\end{array} \right)ds
+\gamma\left( \left( \begin{array}{c}
                   \BOne \\
                   \rho w
\end{array} \right)
+\int_0^t \left( \begin{array}{c}
                   u(s) \\
                   \rho v(s)
\end{array} \right) ds\right)  \, ,
\end{equation}
and solve it for $(\zeta, \dot{\zeta})$ for fixed $t$, 
with $(\zeta(0),\dot{\zeta}(0))=(\BOne, w)$. We obtain the
mapping 
\begin{equation} \label{nohm} 
\mc{M}_F: (\eta(t), \dot{\eta}(t)) \rightarrow (\zeta(t),\dot{\zeta}(t))
\end{equation}
Its fixed point, over the time interval where it is defined, gives the
solution curve to the equations of motion of the body under the action
of $F$, and with the said initial conditions.  

Explicitly, at a general $(\eta, \dot{\eta})$, we have that
$$
\begin{array}{rcl}
D_{(\eta , \dot{\eta })}A_{\lambda }\left( \! \!  
                                  \begin{array}{c}
                                       u \\ v
                                  \end{array}
                            \! \! \right)  & =  & 
  A_{ij}^{\alpha \beta}(D\eta) \partial_i \partial_j u^{\beta} +
(\partial_{(\partial _{k} \eta ^{\gamma })}   
 A_{ij}^{\alpha \beta } )(D\eta) 
(\partial _{i}\partial _{j}\eta ^{\beta })(\partial _{k}u^{
\gamma }) \\ & & +J(\eta)({\rm div}_{\eta}u)\nabla_{\eta}(q+\lambda J(\eta))
+J(\eta) \left[ \bar{u}\cdot \nabla ,\nabla \right]
        _{\eta }(q+ \lambda J(\eta)) \vspace{1mm}  \\ 
& & + J(\eta) \nabla_{\eta }(h+\lambda J(\eta){\rm div}_{\eta}u ) 
\, ,  
\end{array}
$$
where $h$ is defined as the solution to the boundary value problem 
$$
\begin{array}{rcl}
L(\eta )h  & =  &  -{\rm div}_{\eta }(    
  A_{ij}^{\alpha \beta}(D\eta) \partial_i \partial_j u^{\beta} +
(\partial_{(\partial _{k} \eta ^{\gamma })}   
 A_{ij}^{\alpha \beta } )(D\eta) 
(\partial _{i}\partial _{j}\eta ^{\beta })(\partial _{k}u^{
\gamma }) ) - \left[\bar{u}\cdot \nabla ,{\rm div} \right]_{\eta }
{\rm Div}\, W'(D \eta) \vspace{1mm}\\ & & 
- \left[\bar{u}\cdot \nabla ,{\rm div} \right]_{\eta }J(\eta)\nabla_{\eta}q 
- {\rm div}_{\eta}( J(\eta){\rm div}_{\eta} u \nabla_{\eta}q)  
-  {\rm div}_{\eta}J(\eta) \left[\bar{u}\cdot \nabla ,{\rm div} \right]_{\eta } 
q  \vspace{1mm} \\ 
    & & + 2\rho {\rm trace}(-D\bar{u}(D\eta ^{-1}D\dot{\eta })^{2}+
(D \eta )^{-1}D\bar{v}(D\eta )^{-1}D\dot{\eta }) \; , \vspace{1mm} \\
h\mid _{\partial \Omega }  & = & 
{\displaystyle -\frac{1}{J^{b}(\eta)}\left(  A_{ij}^{\alpha \beta}(D\eta)\partial_j u^{\beta}
N^i \nu^\alpha \circ \eta  +
q J^b (\eta) \left[ \det{\left(
 \begin{array}{c}
\partial_{T_{1}} u \\ \partial_{T_{2}}\eta 
 \end{array} \right)}+
 \det{\left( \begin{array}{c} 
   \partial_{T_{1}}\eta \\ \partial_{T_{2}} u 
    \end{array} \right)}
   \right] \right)
 } \, , 
\end{array} 
$$
and where $q$ solves the boundary value problem {\rm (\ref{newq})}.
Here $\bar{u}=u\circ \eta ^{-1}$, $\bar{v}=v\circ \eta ^{-1}$,  and
$\{ T_1, T_2\}$ is an orthonormal frame of the boundary. The boundary
condition for $h$ arises by expressing the boundary condition in (\ref{newq}) as
$$
W'(D\eta) N +qJ^b(\eta) \nu \circ \eta =0 \, , 
$$
and showing that the linearization of $J^b(\eta) \nu \circ \eta$ is 
$J^{b}(\eta)$ times the bracketed sum of determinants in the expression above
for $h\mid_{\partial \Omega}$. Notice that if $\eta =\BOne$, this term is
${\rm div}^{b}(u):=\<\partial_{T_{1}}u,T_{1}\>+
\<\partial_{T_{2}}u,T_{2}\>$.

At $(\eta,\dot{\eta})=(\BOne,0)$, the linearized equation has a simple
expression.  
The linearization of the boundary conditions in (\ref{new2}) at 
$\eta =\BOne$ yield
\begin{equation} \label{nblg}
\begin{array}{c}
{\rm div}\, u \mid_{\partial \Omega} =0 \, ,  \vspace{1mm} \\
(\< \partial _{T}u, W'(\BOne)N \> + \< T,\partial_s(W'(D\eta(s))N)\mid_{s=0} \>
 )\mid_{\partial \Omega} 
= 0\, ,  
\end{array} 
\end{equation}
and so, ${\rm div}^{b}(u)= -\< \partial_N u , N\>$. Then, by evaluating 
$D_{(\BOne,0)}A_{\lambda}$, 
the system (\ref{lah}) reduces to 
\begin{equation}\label{vp}
\frac{d}{dt}\left( \begin{array}{c}
                   u \\ \rho v
                   \end{array}\right) =
      \left(
           \begin{array}{c}
           v \\ A_{ij}^{\cdot \,
\beta}(\BOne) \partial_i \partial_j u^{\beta}+ \nabla (h + \lambda 
{\rm div}\,  u) 
           \end{array}
      \right) +
\left( \begin{array}{c}
                   0 \\ F
       \end{array}\right) \, , 
\end{equation}
where $h$ solves the boundary value problem 
\begin{equation} \label{hah}
\begin{array}{ccl}
\Delta h  & =  &  - {\rm div}A_{ij}^{\alpha \beta}(\BOne) \partial_i 
\partial_j u^{\beta}     \, , \\
h\mid _{\partial \Omega }  & = &
{\displaystyle - A_{ij}^{\alpha \beta}(\BOne )\partial_j u^{\beta}
N^i N^\alpha - \< W'(\BOne)N,N\>\< \partial_N u, N\>  \, . }
\end{array}
\end{equation}

We write the solution of (\ref{vp}) with trivial initial condition as 
$(u(t),v(t))=S_F(\BOne, 0)$. We 
then have defined a right side for 
system (\ref{imp}), whose solution as a path in $t$ we
write as $(\zeta, \dot{\zeta})=S_{G_\lambda}(u(t),v(t))$. The pair
$(\zeta, \dot{\zeta})= (S_{G_\lambda} \circ S_F)(\BOne, 0)=
\mc{M}_F (\BOne, 0)$ is the first orbit point of the 
mapping $\mc{M}_F$ for a body initially stationary at $(\eta(t),\dot{\eta}(0))=
(\BOne, 0)$, and moving subject to the action of the external 
pressure force $F$. 

We use here a numerical evaluation of $(u(t),v(t))=S_F(\BOne, 0)$ 
to approximate the actual nonlinear motion of this  
$\Omega$. This improves our work in \cite{sim}, where the said motion was 
approximated 
by the solution of the linearization of (\ref{eqs}) itself on the submanifold 
defined by (\ref{eq1}) \cite[system (22)]{sim}, and where the actual  
motion needed to be small enough so that it could be approximated well in this 
manner. The changes now, and when numerically
possible, the use of the first iterate $\mc{M}_F(\BOne,0)=S_{G_\lambda}
(S_F(\BOne, 0))$ itself to approximate the motion of $\Omega$, widens 
significantly the range where the approximation is reasonably accurate.

The methods here and in \cite{sim} yield compatible results in the
regime where they are both applicable, but the differences impose 
significant changes when it comes to the numerical evaluation of the 
system (\ref{vp}) now involved. The finding of its numerical solution 
requires a discretizing space of richer structure 
than the one used in treating \cite[system (22)]{sim}. Since
$\lambda \neq 0$, the $u(t)$ that we seek now is not necessarily 
a divergence-free vector field, 
although it is close to one given the initial conditions, and tends to 
be driven even closer to one by the equation it satisfies, as ultimately, the
actual solution of the equations of motion has a divergence-free velocity.  
Our work here is thus harder than that in \cite{sim}.

This extra effort in our work now is justified by the better accuracy of the 
approximation to the actual motion that we obtain, and by the fact that if the  
procedure were to be iterated (with the subsequent linearizations 
carried out at the previously found $(\zeta,\dot{\zeta})$), we would produce
a sequence that converges to the solution of the nonlinear elastic 
motion on some time interval, a possibility not available
when using the numerical scheme in \cite{sim}.
Going a bit further, if we were to 
take Cauchy data for (\ref{vp}) that is compatible 
with a nonzero divergence-free initial velocity $\dot{\eta}(0)=w$, 
the solution $(u,v)=S_F(\BOne, w)$, and  
the ensuing orbit point $(\zeta,\dot{\zeta})=
(S_{G_\lambda}\circ S_F)(\BOne, w)$
would serve to describe the motion of an initially moving $\Omega$ with
velocity $w$ under 
the influence of the external pressure force $F$.
   
\section{Hooke bodies: Orthotropic materials}
\label{s3}
We let $\mc{S}^2$ denote the space of symmetric 2-tensors on $\Omega$, and
$\sigma$ and $e$ be the stress and strain tensors respectively. 
Since we assume conservation of momentum, the tensor $\sigma $ is symmetric.
Its components have the dimension of force per unit 
area, or pressure. The tensor $e$ is symmetric; if $u$ is the
displacement $\eta (t)(x)=x+u(t,x)$, and we use the 
Euclidean metric in $\mb{R}^3$, we have that
\begin{equation} \label{eq13}
e_{ij}=\frac{1}{2}(\nabla u +\nabla u^T - \nabla u^T \nabla u)_{ij} \, ,
\end{equation}
and modulo quadratic errors, $e$ coincides with the symmetrized vector of 
covariant derivatives of $u$. We often equate the two; the latter notion is 
usually called the infinitesimal strain. The components of $e$ are 
dimensionless.

The body $\Omega$ is said to be of Hooke type if there exists a tensor
$W \in {\rm End}(\mc{S}^2)$,   
$$
W: C^{\infty}(\Omega; \mc{S}^2) \rightarrow C^{\infty}(\Omega; \mc{S}^2)
$$
such that  $\sigma = W e$, 
and whose stored energy function is given by 
\begin{equation} \label{enc}
W(D \eta)=\frac{1}{2}\< \sigma, e\>=\frac{1}{2}\< We,e\>\, . 
\end{equation}
This tensor $W$ is called the tensor of elastic constants, 
or moduli, of the material. 

In components $W=(W^{ijkl})$, we have that
\begin{equation} \label{dec}
\sigma^{ij}=W^{ijkl} e_{kl} \, ,  
\end{equation}
and
\begin{equation} \label{sef} 
W(D \eta)=\frac{1}{2}\< \sigma, e\>=\frac{1}{2}W^{ijkl} e_{kl} e_{ij}\, , 
\end{equation}
with the symmetries $W^{ijkl}=W^{jikl}=W^{ijlk}$.
Coercivity of $W$ imposes the additional symmetry $W^{ijkl}=W^{klij}$, 
yielding a total of $21$ degrees of freedom for
$W$. Explicitly, we have 
\begin{equation} \label{ten}
\left( \begin{array}{c}
\sigma_{11} \\
\sigma_{22} \\
\sigma_{33} \\
\sigma_{23} \\
\sigma_{31} \\
\sigma_{12} 
\end{array}
\right)
= 
\left( \begin{array}{cccccc}
W_{11}^{\phantom{1}11} &
W_{11}^{\phantom{1}22} &
W_{11}^{\phantom{1}33} &
W_{11}^{\phantom{1}23} &
W_{11}^{\phantom{1}31} &
W_{11}^{\phantom{1}12} \\
W_{22}^{\phantom{1}11} &
W_{22}^{\phantom{1}22} &
W_{22}^{\phantom{1}33} &
W_{22}^{\phantom{1}23} &
W_{22}^{\phantom{1}31} &
W_{22}^{\phantom{1}12} \\
W_{33}^{\phantom{1}11} &
W_{33}^{\phantom{1}22} &
W_{33}^{\phantom{1}33} &
W_{33}^{\phantom{1}23} &
W_{33}^{\phantom{1}31} &
W_{33}^{\phantom{1}12} \\
W_{23}^{\phantom{1}11} &
W_{23}^{\phantom{1}22} &
W_{23}^{\phantom{1}33} &
W_{23}^{\phantom{1}23} &
W_{23}^{\phantom{1}31} &
W_{23}^{\phantom{1}12} \\
W_{31}^{\phantom{1}11} &
W_{31}^{\phantom{1}22} &
W_{31}^{\phantom{1}33} &
W_{31}^{\phantom{1}23} &
W_{31}^{\phantom{1}31} &
W_{31}^{\phantom{1}12} \\
W_{12}^{\phantom{1}11} &
W_{12}^{\phantom{1}22} &
W_{12}^{\phantom{1}33} &
W_{12}^{\phantom{1}23} &
W_{12}^{\phantom{1}31} &
W_{12}^{\phantom{1}12} 
\end{array}
\right)
\left( \begin{array}{c}
e_{11} \\
e_{22} \\
e_{33} \\
e_{23} \\
e_{31} \\
e_{12} 
\end{array}
\right) \, .
\end{equation}

For bodies of Hooke type, we have that   
\begin{equation} \label{bs}
\begin{array}{rcl}
W'(\BOne)^{i}_{\mbox{}\hspace{1mm}\alpha} & = & {\displaystyle 
\frac{\partial W}{\partial( \partial_i \eta^\alpha)}\mid_{\eta =\BOne}=
W^{i \hspace{1mm} j}_{\mbox{}\hspace{1mm}\alpha\hspace{1mm}j}} \, , \vspace{1mm}
\\
\partial_s(W'(D\eta(s))N)\mid_{s=0} & = &  
W^{\alpha \hspace{1mm} j}_{\mbox{}\hspace{1.5mm}i \hspace{1.5mm}\beta} 
\partial_j 
u^\beta N^i=\sigma(\nabla u)^{\alpha}_{\mbox{} \hspace{1mm} i}N^i=
\sigma(\nabla u)N \, , \vspace{1mm} \\
A_{ij}^{\alpha \beta}( \BOne ) & = &
W^{i\hspace{1.5mm}j}_{\phantom{i}\alpha \hspace{1.5mm} \beta}\, . 
\end{array}
\end{equation}

Orthotropic materials are bodies of Hooke type that posses three mutually 
orthogonal planes of symmetries at each point, with three corresponding 
orthogonal axes, and so have unchanging elastic coefficients  
under rotations of $180^{\circ}$ about any of these axes. Consequently, 
the tensor $W$ of these bodies has only $9$ 
degrees of freedom, and its expression (\ref{ten}) 
relative to these preferred axes reduces to 
\begin{equation}\label{ten2}
W=\left( \begin{array}{cccccc}
W_{11}^{\phantom{1}11} &
W_{11}^{\phantom{1}22} &
W_{11}^{\phantom{1}33} &
0 & 0 & 0  \\
W_{22}^{\phantom{1}11} &
W_{22}^{\phantom{1}22} &
W_{22}^{\phantom{1}33} &
0 & 0 & 0  \\
W_{33}^{\phantom{1}11} &
W_{33}^{\phantom{1}22} &
W_{33}^{\phantom{1}33} &
0 & 0 & 0  \\
0  & 0  & 0 &
W_{23}^{\phantom{1}23} &
0 & 0 \\
0  & 0  & 0 & 0  &
W_{31}^{\phantom{1}31} &
0 \\
0 & 0 & 0 & 0 & 0 &
W_{12}^{\phantom{1}12} 
\end{array}
\right) \, .
\end{equation} 

The elastic constants of an orthotropic Hooke body are parametrized by   
the three moduli of elasticity, the six Poisson ratios, and the three
moduli of rigidity, or shear modulus, determined by the axes of symmetry.
The moduli of elasticity $e_1, e_2, e_3$, and moduli of rigidity $g_1, g_2,
g_3$, have dimensions of force per unit area, while the 
Poisson ratios $\mu_{12}, \mu_{21}, \mu_{13},\mu_{31}, \mu_{23},
\mu_{32}$ are dimensionless. The compatibility
relations
\begin{equation} \label{rel} 
\frac{\mu_{ij}}{e_i}=\frac{\mu_{ji}}{e_j}\, , \quad \text{$i\neq j$, 
$1\leq i,j \leq 3$}\, , 
\end{equation}
leaves a total of $9$ independent parameters. 

Wood is considered as a typical example of orthotropic material  
since it has unique, and somewhat independent, mechanical properties along
three mutually perpendicular directions: The longitudinal axis $z$ that is
parallel to the fiber grains; the radial axis $r$ that is normal to the growth
rings; and the tangential axis $\theta$ to the growth rings.
If we place the vertical axis of a tree trunk along the $z$-direction, 
the axes of symmetry of its wood coincide with the ordered 
cylindrical coordinates $(r,\theta, z)$. 
Then the elastic behaviour of wood 
is described by the elasticity moduli $e_r, e_\theta, e_z$, the rigidity 
moduli $g_{z r}, g_{z\theta}, g_{r \theta}$, and six Poisson ratios $\mu_{r 
\theta},\mu_{\theta r}, \mu_{r z}, \mu_{z r}, \mu_{\theta z}, \mu_{z \theta}$.
These constants satisfy the three relations (\ref{rel}), and their
relation to the components of the moduli tensor is explicitly given by
\begin{equation} \label{decI}
U=W^{-1}= \left( \begin{array}{cccccc}
{\displaystyle \frac{1}{e_r} } &
{\displaystyle -\frac{\mu_{\theta  r}}{e_\theta}}  &
{\displaystyle -\frac{\mu_{z r}}{e_z}}  &
0 & 0 & 0  \vspace{1mm} \\
{\displaystyle -\frac{\mu_{r \theta }}{e_r}} &
{\displaystyle \frac{1}{e_\theta}}  &
{\displaystyle -\frac{\mu_{z \theta}}{e_z}}  &
0 & 0 & 0  \vspace{1mm} \\
{\displaystyle -\frac{\mu_{r z}}{e_r}} &
{\displaystyle -\frac{\mu_{\theta z}}{e_\theta }} &
{\displaystyle \frac{1}{e_z}} &
0 & 0 & 0  \vspace{1mm} \\
0  & 0  & 0 &
{\displaystyle \frac{1}{g_{\theta z}} } &
0 & 0 \vspace{1mm} \\
0  & 0  & 0 & 0  &
{\displaystyle \frac{1}{g_{rz}}} &
0 \vspace{1mm} \\
0 & 0 & 0 & 0 & 0 &
{\displaystyle \frac{1}{g_{r\theta} }} 
\end{array}
\right) 
\, .
\end{equation}
where
\begin{equation} \label{decin}
e^{ij}= U^{ijkl}\sigma_{kl} \, , 
\end{equation}
is the inverse of the tensorial relation (\ref{dec}). 

For lack of better choices, in all of our simulations below, 
we use the material constants for the Engelmann Spruce extracted from 
data provided in \cite{gwk} to obtain the moduli tensor of constant of the 
bodies. These values are  shown in Tables 1 \& 2 below.
We have used these same constants previously in 
\cite{sim}, so we may now draw comparisons of the results, and judge the 
improvement obtained.

\begin{table}
\begin{tabular}{|l|c|c|c|c|c|c|} \hline \hline
  & $e_z$ in MPa & $e_{\theta}/e_z$ & $e_{r}/e_z$ & $g_{z r}/e_z$ &  
 $g_{z\theta}/e_z$ & $g_{r \theta}/e_z$  \\
\hline \hline
Spruce, Engelmann & 9,790 & 0.059  & 0.128  & 0.124  & 0.120  & 0.010  \\
\hline \hline 
\end{tabular}
\medskip 

\caption{Ratios of elasticity to rigidity moduli for 
Engelmann spruce.}
\end{table}

\begin{table}
\begin{tabular}{|l|c|c|c|c|c|c|} \hline \hline
  & $\mu_{z r}$ & $\mu_{z\theta}$ & $\mu_{r \theta}$ & $\mu_{\theta r}$ & $\mu_{r z}$ & 
$\mu_{\theta z}$  \\
\hline \hline
Spruce, Engelmann & 0.422  & 0.462 & 0.530 & 0.255 & 0.083 & 0.058  \\
\hline \hline 
\end{tabular}
\medskip

\caption{Poisson ratios for Engelmann spruce.}
\end{table}
\medskip

The bodies we analyze are thin plates made of spruce with these elastic 
constants, the thinness condition making it natural to assume that the 
components of the elastic tensor are the same when expressed in cylindrical 
or Cartesian coordinates, which we take here as a fact 
(see \cite[Remark 5]{sim}). 
In addition, we shall assume that these components are constant throughout 
the body. 
The values in Tables 1 and 2 reflect a failing 
condition (\ref{rel}), so we take the average of the computed values of
$\mu_{ij}/e_i$ and $\mu_{ji}/e_j$ as the value of either one of 
these quantities in our calculations. Thus, in all of our simulations,
the diagonal blocks of the tensor $W$ in (\ref{ten2}) are 
\begin{equation} \label{ten3}
\begin{array}{c}
W_{3\times 3}=10^7 \left( \begin{array}{ccc}
157.198269069862 & 44.1920517114940 &  116.065341927474  \\ 
44.1920517114940 & 72.0200103705017 & 75.6887031695923  \\ 
116.065341927474 & 75.6887031695923 & 1095.80735919001  
\end{array}\right) \, , \vspace{1mm} \\
D_{3\times 3} =10^7 \left( \begin{array}{ccc}
117.480 & 0 & 0 \\ 
0 & 121.396 & 0 \\ 
0 & 0 & 9.790 
\end{array}\right) \, ,
\end{array}
\end{equation}
respectively, the unit of measurement Pa.  
As we work in a Cartesian orthonormal frame, we can raise or lower indices 
in tensors with abandon. We take $\rho = 360 \; {\rm kg}/{\rm m}^3$ for the 
density parameter.  Notice that the eigenvalues of $W_{3\times 3}$ are 
$$ 
10^7\{ 156.292790395160, 1116.15681336097, 52.5760348742398 \} \, , 
$$
so the stored energy function of any of our bodies is coercive.

\section{Smooth triangulations and induced discretizing spaces}
\label{s4}
The geometric nature of the unknowns in the systems (\ref{vp}) and (\ref{imp})
 makes it natural to cast their solutions using 
the Whitney forms of an oriented triangulation of the body, or their
metric duals. We recall these notions briefly.
For general definitions, and properties of the Whitney forms, we refer the
reader to \cite{do}; some additional motivation behind our choices, or reasons
for making them, may be found in \cite[\S 4.1]{sim}.

We consider a connected Riemannian $n$-manifold with boundary 
$(M^n, g)$. In our work, $n=3$, $M^3$ is embedded in $\mb{R}^3$, and $g$ is
the metric induced on it by the Euclidean metric in the ambient space; these
$M^3$s are oriented, with their orientation compatible with that of 
$\mb{R}^3$.

We let $K$ be a finite smooth oriented triangulation of $M$.
We identify the polytope of $K$ with $M$, and fix some ordering of the 
vertices of $K$. We denote by $C^q(K)$ the space of simplicial {\it oriented} 
$q$-cochains, and by $L^2(M;\Lambda^q T^*M)$ the space of $L^2$ $q$-forms on 
$M$. If $P \in C^{\infty}(M; E) \rightarrow C^{\infty}(M;\tilde{E})$ is a 
linear operator over $M$ mapping sections of $E$ to sections of $\tilde{E}$,
we denote by $L^2_P(M; E)$ the subspace of $L^2(M;E)$ sections that are mapped
by $P$ into $L^2(M;\tilde{E})$, provided with the graph norm.     
  
The barycentric subdivision $K'$ of a (not necessarily oriented) triangulation
$K$ is a simplicial complex that is naturally oriented, its  
vertices ordered by decreasing dimension of the 
simplices of the triangulation of which they are the barycenters. This 
ordering induces a linear ordering of the vertices of each simplex of $K'$
\cite{mu}.  

If $K'$ is the the barycentric subdivision 
of the smooth triangulation $K$ of the Riemannian $3$-manifold $(M^3,g)$, we  
denote by $K^{(j)}$ its $j$th skeleton,  
and by $\mc{V}'$, $\mc{E}'$, $\mc{F}'$, and $\mc{T}'$, the set of vertices, 
edges, faces, and tetrahedrons  
of $K'$, respectively. The interior edges in $K^{(1)}$ 
are denoted by $\mc{E}'_{\circ}$, while the boundary edges 
are denoted by $\mc{E}'_{\partial }$. The interior faces in $K^{(2)}$  
are denoted by $\mc{F}'_\circ$, while the boundary faces
are denoted by $\mc{F}'_\partial$. 
If necessary, the analogous concepts for the triangulation $K$ itself will 
be denoted similarly but without the $'$s. Notice that $K^{(1)}$ is the 
oriented graph $G(\mc{V}', \mc{E}')$. 
The cardinality of a set $S$ is denoted by $|S|$.  

Any oriented triangulation of $(M^n,g)$ has associated with it 
the set of piecewise linear Whitney forms \cite{whi}, and their  
corresponding metric duals. When $n=3$, the metric duals of the
said forms are functions in degree zero and three, and vector fields in 
degree one and two, respectively. They all play roles in our work.
We use the barycentric subdivision $K'$ of the triangulation $K$, and 
describe these Whitney forms, and their metric duals, in that particular case. 
 
If $p \in \mc{V}'$, we let $x_p$ be the $p$-th barycentric coordinate function
in $K'$. This is the Whitney form of degree zero associated to the vertex $p$.
The collection $\cup_p \{ x_p \}_{p \in \mc{V}'}$ is a
partition of unity of the polytope of $K'$ ``subordinated'' to the open 
cover $\{ {\rm St}\, p\}$. Notice that
in the weak sense, $\nabla^g x_p$ is a well-defined $L^2$-vector field. 
We define the space  
\begin{equation} \label{fun}
{\rm Fun}(K') = {\rm span} \{ x_p \}_{p \in \mc{V}'} \, .
\end{equation}
We have that ${\rm dim}\,{\rm Fun}(K')=|\mc{V}'|$. 
 
If $e=[p_0,p_1] \in \mc{E}'$, we consider the piecewise continuous $1$-form
$w_e= x_{p_0} d x_{p_1}- x_{p_1}dx_{p_0}$, the Whitney form associated to $e$,
and its metric dual vector field
\begin{equation} \label{f1e}
W_e= x_{p_0} \nabla^g x_{p_1}- x_{p_1}\nabla^g x_{p_0} \, .
\end{equation}
It is an element of the space of $L^2$ forms, or vector fields, with integrable 
squared norm. It has compact support in 
$\overline{\rm St}\, p_0 \cap \overline{\rm St}\, p_1=\overline{{\rm St}}\, e$,
and in the weak sense, ${\rm curl} \, W_e$ is a well-defined $L^2$-vector 
field. Further, for any oriented edge $e'$ of $K'$, 
$$ 
e^*(e')=\< w_e, e'\> =\int_{e'} W_e = \left\{ \begin{array}{ll}
                      1 & \text{if $e'= e$}\, , \\
                      0 & \text{otherwise}\, ,
                      \end{array}
              \right.    
$$ 
and if $b_{pe}$ is the incidence number of the vertex $p$ 
and edge $e$ in the graph $G(\mc{V}',\mc{E}')=K^{(1)}$,  
we have that
\begin{equation} \label{grf}
U_p:=\nabla^g x_p =  
\sum_{q\in \mc{V}'} (x_q \nabla^g x_p -
x_p\nabla^g x_q)=(\sum_{q\in \mc{V}'} x_{q})\nabla^g x_p -
x_p \nabla^g ( \sum_{q\in \mc{V}'} x_q )   = \sum_{e\in \mc{E}'} b_{pe}W_e \, . 
\end{equation}
(This last identity implies also that ${\rm curl}\,  U_p$ is well-defined, 
and identically zero.)

We denote the spaces spanned by the $W_e$s in (\ref{f1e}) and by the $U_p$s in
(\ref{grf}) as
\begin{equation} \label{grad}
\begin{array}{rcl}
{\rm Grad}(K') & = & {\rm span }\{ W_e\}_{e\in \mc{E}'} \, , \\
{\rm Grad}_0(K') & = & {\rm span}\{ U_p \}_{p\in \mc{V}'} \, .
\end{array}
\end{equation}
Their dimensions are  ${\rm dim}\, {\rm Grad}(K')=| \mc{E}'|$ and
${\rm dim}\, {\rm Grad}_0(K')=| \mc{V}'|-1$, respectively.

The Whitney forms of degree two are constructed from the faces of $K'$. For if 
$f=[p_0,p_1,p_2] \in \mc{F}'$, we associate with it the $L^2$ form 
$w_f= 2(x_{p_0}dx_{p_1} \wedge d x_{p_2}+x_{p_1}dx_{p_2} 
\wedge d x_{p_0} +x_{p_2}d x_{p_0} \wedge d x_{p_1})$. Its Hodge
$*$ is a degree one form whose metric dual is 
the vector field
\begin{equation} \label{f2f}
W_f= 2(x_{p_0}\nabla^g x_{p_1} \times \nabla^g x_{p_2}+x_{p_1}\nabla^g x_{p_2} 
\times \nabla^g x_{p_0} +x_{p_2}\nabla^g x_{p_0} \times \nabla^g x_{p_1}) \, , 
\end{equation}
where $\nabla^g x_{p_i} \times \nabla^g x_{p_j}$ is the cross product
of the $L^2$-gradients $\nabla^g x_{p_i}$ and $\nabla^g x_{p_j}$, 
respectively. 

The flux of $W_f$ through a face $f'$ in $\mc{F}'$ is 
well-defined, and given by  
$$ 
f^*(f')=\int_{f'} W_f \cdot n_{f'}  = \left\{ \begin{array}{ll}
                      1 & \text{if $f' = f$}\, , \\
                      0 & \text{otherwise}\, .
                      \end{array}
              \right.
$$ 

The family of vector fields $\{ W_f \}_{\text{$f\in \mc{F}'$}}$ 
is linearly independent, and its span
contains the image under {\rm curl} of ${\rm Grad}(K')$.
Indeed, if $e=[p_0,p_1]$ is an edge in $\mc{F}'$, 
and $b_{ef}$ is now the incidence number of the edge 
$e$ on the face $f$, we have that 
\begin{equation} \label{vf}
U_e:= {\rm curl}\, W_e= 2 \nabla^g x_{p_0}
\times \nabla^g x_{p_1} = \sum_{f\in \mc{F}'} b_{ef}W_f \, .
\end{equation}
(It follows from this identity that the weak divergence of $U_e$ in $L^2$ is 
well-defined, and identically zero.)
Notice in addition that if $b_{ft}$ is the incidence number of the face $f$ in
the tetrahedron $t$, we have that
$$
{\rm div} \, (\sum_{f\in \mc{F}'}c_f W_f) =\sum_{t\in \mc{T}'}
(\sum_{f \in \mc{F}'} c_f b_{ft}) t^{*}\, ,
$$
where $t^*$ is the characteristic function ordered $3$-cochain determined
by the tetrahedron (or $3$-simplex) $t$.  Thus, 
$\sum_{f\in \mc{F}'}c_f W_f$ is divergence-free if, and only if, the
weighted sum $\sum_{f\in \mc{F}'}c_f b_{ft}$ over the four faces of each
tetrahedron $t$ in $K'$ is identically zero.

We denote the spaces spanned by the $W_f$s in (\ref{f2f}) and by the $U_e$s in
(\ref{vf}) as
\begin{equation} \label{div}
\begin{array}{rcl}
{\rm Div}(K') & = & {\rm span }\{ W_f\}_{f\in \mc{F}'} \, , \\
{\rm Div}_0(K') & = & {\rm span}\{ U_e \}_{\text{$e\in \mc{E}'$}} \, .
\end{array}
\end{equation}
Their dimensions are  ${\rm dim}\, {\rm Div}(K')=| \mc{F}'|$ and
${\rm dim} \, {\rm Div}_0(K')= | \mc{E}'|-(|\mc{V}|'-1)$, respectively.

Finally, the Whitney form of a tetrahedron $t=[p_0,p_1,p_2,p_3] \in K'$ is
defined by $w_t=6(x_{p_0}dx_{p_1}\wedge dx_{p_2} \wedge dx_{p_3} 
-x_{p_1}dx_{p_0}\wedge dx_{p_2} \wedge dx_{p_3} 
+x_{p_2}dx_{p_0}\wedge dx_{p_1} \wedge dx_{p_3} 
-x_{p_3}dx_{p_0}\wedge dx_{p_1} \wedge dx_{p_2})$. It is a piecewise linear 
three form supported on $t$, and since the restriction of 
$x_{p_0}+x_{p_1}+x_{p_2}+x_{p_3}$ to the polytope $|t|$ is equal to the
constant function one, we can identified it with the $3$-form with 
$6 d\mu_t$, $d\mu_t$ the natural volume form of the $3$-simplex $t$. By 
the Hodge star operator, we then see that 
$w_t$ corresponds to the $3$-cochain locally constant function $W_t=6 t^{*}$.
The set $\{ \frac{1}{6}W_t\}_{t \in \mc{T}'}$ is linearly independent, and the
space 
\begin{equation}\label{char}
{\rm Char}(K')= {\rm span} \left\{ \frac{1}{6} W_t\right\}_{t \in \mc{T}'}
\end{equation}
that it spans is a subspace of $L^2$ of dimension $| \mc{T}'|$. It constitute
a partition of unity of the polytope of $K'$ through locally
constant functions ``subordinated'' to the
covering $\{ \overline{t}\}_{t \in \mc{T}'}$.

The space ${\rm Fun}\, (K')$ is used in the usual manner to discretize
scalar valued functions in $L^2(M)$ as a combination of locally 
supported continuous terms, 
its basis elements encoding the combinatorial property of all adding to the
constant function $1$. The space ${\rm Grad}(K')$ is the
natural choice for discretizing gradient vector fields in $L^2(M;TM)$. Its
subspace ${\rm Grad}_0(K')$ is spanned by piecewise constant vector fields whose
basis elements are true gradients. But they yield trivial results when they  
are acted on by differential operators of nonzero order that annihilate the 
constants, thus making it a necessity to enlarge the view, and consider
${\rm Grad}(K')$ instead. Similarly, the 
space ${\rm Div}(K')$ is the natural choice when discretizing 
divergence-free vector fields in $L^2(M;TM)$, its 
subspace ${\rm Div}_0(K')$ consisting of elements that though divergence-free 
per se, are piecewise constant and so acted on by differential operators of 
nonzero order in a trivial manner. Finally, 
${\rm Char}(K')$ is the natural choice as discretizing space for 
the divergence of vector fields as a combination of locally constant terms, 
since ${\rm div}({\rm Div}(K')) \subset {\rm Char}(K')$.
(Notice that by (\ref{f1e}), we have that
${\rm div}\, {\rm Grad}(K')=0$ in the $L^2$-sense.)

In any of these cases, the consistency between the geometric property of the 
discretization of the scalar or vector fields, and the choices of space where 
it is carried out, is encoded in the adjacency matrices of
the triangulation in use, which in turn is a reflection of the fact that 
the homology of the complex
$$
0 \rightarrow C^{\infty}(M) \stackrel{{\rm grad}}{\longrightarrow}
C^\infty(M,TM) \stackrel{{\rm curl}}{\longrightarrow} C^\infty(M,TM)
\stackrel{{\rm div}}{\longrightarrow}  C^\infty(M) \rightarrow 0
$$
equals the cohomology of $M$, and can be computed from the cohomology of
its discretized $L^2(K')$-version
$$
0 \rightarrow F(K') \stackrel{{\rm grad}}{\longrightarrow}
{\rm Grad}(K') \stackrel{{\rm curl}}{\longrightarrow} {\rm Div}(K') 
\stackrel{{\rm div}}{\longrightarrow}  {\rm Char}(K')\rightarrow 0\, .
$$
For the bodies of interest to us, the polytope of $K$ (and, consequently, of 
$K'$) is contractible to a point, or to the wedge of two circles. Thus, 
the kernel 
of the divergence operator on ${\rm Div}(K')$ coincides with ${\rm Div}_0(K')$;
in general, though, this is true only modulo a finite dimensional space whose
dimension is the rank of the second cohomology group of $K'$.
In degree one, the kernel of the curl operator on ${\rm Grad}(K')$ agrees with
${\rm Grad}_0(K')$ if, and only if, the first cohomology of $M$ is trivial;
otherwise, the equality holds modulo a finite dimensional space whose
dimension is the rank of the first cohomology group of $K'$. The cohomology
groups in degrees zero and three have rank one, the cycle in both cases
being the constant function $1$ expressed as $1=\sum_{p\in \mc{V}'} x_p$,
and $1 =\sum_{t\in \mc{T}'}\frac{1}{6} W_t$, respectively.   

Although we ultimately seek solutions to the equations of motion (\ref{eqs}),
and these are given by curves of diffeomorphisms whose tangent vectors are
divergence free fields, the intermediate steps in solving (\ref{mod}) to get
to these solutions produce vector fields that do not have this property. The 
linearized equation (\ref{vp}) that we solve here, and its analogue in 
\cite{sim}, contrast in that respect. 
As we solve numerically a weak version of (\ref{vp}) in the Sobolev space 
$H^1(\Omega;\mb{R}^3)$, or in $C^{1,\alpha}(\Omega;\mb{R}^3)$, 
the space of choice for discretizing the sought after solution is
\begin{equation} \label{l2s}
L^2_1(K';\mb{R}^3):={\rm Grad}(K') \oplus {\rm Div}(K') \, ,
\end{equation}
with the discretization expressed in the natural Whitney basis elements 
defining the summands $\mc{B}_{L^2_1}=\{ W_s \}_{W_s \in L^2_1}:=
\{ W_e, W_f \}_{(e,f) \in \mc{E}'\times \mc{F}'}$. 
This space is really the basic set-up for the proof of the duality
$H^1(M;\mb{R}) \cong H_2(M;\mb{R})$, one step away from the sum of 
simplicial and dual block decompositions of $K'$ from where this proof 
departs, and corresponds to the $L^2$ decomposition of a vector 
field into a gradient plus a divergence-free component. (Analogously, the space 
$$
{\rm Fun}(K') \oplus {\rm Char}(K') \, ,
$$
in which it would be natural to discretize any scalar valued function 
defined on the polytope $|K'|$, corresponds to the $L^2$ decomposition of 
a function as one in the image of the Laplace operator plus its
projection onto the constants, and in a sense correlated to the one above, 
it is the basic set-up for the proof of the duality 
$H^0(M;\mb{R}) \cong H_1(M;\mb{R})$.)

We emphasize the fact that although in the problems treated here the manifold
$M$ is oriented, all of the spaces defined above do not depend on that,
and it is only the orientation of the simplicial complex $K'$ that 
matters. The latter allows for the fixing of compatible local 
orientations nearby any simplex in the complex $K'$ that, if $M$ 
were to be oriented, would be compatible with this global orientation 
\cite{sidc}; this local orientation is all that is required to carry out 
the Hodge * operation on the Whitney form associated to any simplex.
Once we think about it, this situation is very natural;
it becomes transparent when, for example, we attempt to study the motion of
incompressible perfect fluids, or incompressible elastic bodies in general, on 
nonoriented Riemannian manifolds. The framework developed above 
for the discretizing spaces works verbatim in that context also, oblivious
to this global orientation issue, and merely requiring the choice of a 
$1$-density on $M$ that can be used to define the discrete 
$L^2$-spaces above, and that would have to have been given anyway 
in order to define the Lagrangian (\ref{lag}) that would get the whole 
theory started.
 
We continue our work often relaxing, without mentioning it,  
the smoothness assumption on $M$ to
that of being a smooth manifold with corners. All of the spaces above 
associated  to the triangulation $K$, as well as the $L^2$ spaces that were 
considered, have natural extensions to that context if (some of) 
the differential operators involve in their definition are interpreted weakly.

\section{The algorithms}
\label{s5}
We write the Cauchy problem for the linear system   
(\ref{vp}) as 
\begin{equation}\label{vp2}
\begin{array}{lcl}
{\displaystyle \rho \frac{d^{2}u}{dt^{2}}} & = & A_{ij}^{\cdot \, \beta}(\BOne)
\partial_i \partial_j u^\beta + \nabla h +\lambda \nabla {\rm div}\, u
+F \, , \vspace{1mm} \\
u\mid_{t=0} & = & u_0 \, , \vspace{1mm}  \\
\partial_t u\mid_{t=0} & = & u_1 \, ,
\end{array}
\end{equation}
where, by (\ref{bs}), $A_{ij}^{\alpha \, \beta}(\BOne)=
W^{i\hspace{1.5mm}j}_{\phantom{i}\alpha \hspace{1.5mm} \beta}$, and  
equation (\ref{hah}) for $h$ reduces to 
\begin{equation}\label{hb}
\begin{array}{rcl}
\Delta h & = & -{\rm div}\,
(W^{i\hspace{1.5mm}j}_{\phantom{i}\alpha \hspace{1.5mm} \beta} \partial_i
\partial_j u^\beta )\\
h\mid_{\partial \Omega} & = &  -N^i N^\alpha W^{\phantom{i}\hspace{1.5mm}j}_
{i\alpha \hspace{1.5mm} \beta} \partial_j u^\beta -
W^{\phantom{i}\hspace{1.5mm}j}_{i\alpha \hspace{1.5mm} j}N^i N^{\alpha}
\< \partial_N u, N\> = -\< \sigma(\nabla u)N, N\> -
 \< W'(\BOne)N, N\>\<\partial_N u, N\>\, .
\end{array}
\end{equation}
We analyze numerically its weak solutions in $H^1(\Omega;\mb{R}^3)$.

The reader should notice that the trace map $H^1(\Omega;\mb{R}^3) \ni u 
\rightarrow h\mid_{\partial \Omega} \in H^{-\frac{1}{2}}(\partial \Omega)$ 
for the boundary condition in (\ref{hb}) is
not continuous \cite[Corollary 2.3.5]{sibo}, which forces a careful 
interpretation of the meaning
of the $L^2$ pairing $-\< h, {\rm div}\, u\>$ that arises when 
dualizing the term $\nabla h$ in the right side of equation (\ref{vp2})
\cite[\S 2]{ebsi2}. Theoretically, this is resolved by defining 
the boundary value problem (\ref{hb}) first over the
dense subset of $H^1(\Omega;\mb{R}^3)$ consisting 
of the $u$s in $H^2(\Omega;\mb{R}^3)$ that satisfy the boundary 
conditions (\ref{nblg}) (these $u$s are in $H^1_\Delta(\Omega;\mb{R}^3)$)
\cite[\S 2.8.1]{lima}, \cite{sith},  
 and then extending it by continuity to the whole of
$H^1(\Omega; \mb{R}^3)$ \cite[Proposition 2.3.6]{sibo}. 
The solution operator $H^1(\Omega;\mb{R}^3) \ni u \rightarrow h (u) \in 
L^2(\Omega)$ that results is continuous.

The discretization $U$ of a weak solution $u \in H^1(\Omega;\mb{R}^3)$ is 
carried out  
over the space $L_1^2(K';\mb{R}^3)$ in (\ref{l2s}). We use the basis 
$\mc{B}_{L^2_1}$ of this space given by the  families of vector fields 
$\{ W_e \}_{e\in \mc{E}'}$, and $\{ W_f \}_{f\in \mc{F}'}$, respectively.
By the decomposition 
$\mc{E}'=\mc{E}'_{\circ} + \mc{E}'_{\partial}$ and
$\mc{F}'=\mc{F}'_{\circ} + \mc{F}'_{\partial}$ for the edges and faces of
$K'$, we split $U$ into blocks accordingly,
$$
U = \sum_{e\in \mc{E}'} c_{e} W_{e}+ \sum_{f\in \mc{F}'} c_{f} W_{f} =
 \sum_{e_\circ \in \mc{E}'_{\circ}}
 c_{e_\circ} W_{e_\circ }+ 
\sum_{e_\partial \in \mc{E}'_\partial}  c_{e_{\partial}} W_{e_{\partial}} +
\sum_{f_\circ \in \mc{F}'_{\circ}} c_{f_\circ} W_{f_\circ }+ 
\sum_{f_\partial \in \mc{F}'_{\partial }} c_{f_{\partial}} W_{f_{\partial}} \, ,
$$
which in the spirit of Einstein summation convention, we express   
succinctly as 
\begin{equation} \label{asb}
U= c_{e}W_{e}+ c_{f}W_{f}=c_{e_\circ}W_{e_\circ}+ c_{e_\partial}W_{e_\partial}+ 
 c_{f_\circ}W_{f_\circ}+ c_{f_\partial}W_{f_\partial}\, .
\end{equation}
The vector of coordinates 
$\left( \begin{array}{c}
        c_{e} \\ c_{f}
        \end{array}
\right)
$
is found as the solution of the second order differential equation that 
results from the weak formulation of the equation, after making an appropriate
choice for a discretization $h(U)$ in $L^2(K')$ of the function $h(u)$ that
solves the boundary value problem (\ref{hb}).
    
For convenience, we use an {\it orthonormal} frame to write the 
components $A_{ij}^{\alpha \beta}(\BOne)=
W^{i\hspace{1.5mm}j}_{\phantom{i}\alpha \hspace{1.5mm} \beta}$ of the tensor 
$A$. We assume that this tensor is {\it covariantly constant} over the support 
of each of the basis vectors in $L^2_1(K';\mb{R}^3)$, as indicated earlier.
 
\subsection{The boundary conditions} \label{sbc}
The global condition ${\rm div} \, u \mid_{\partial \Omega}=0$ in (\ref{nblg}) 
is enforced always when deriving the weak version of the equation 
(\ref{vp2}). As for the remaining conditions in (\ref{nblg}), of a local 
nature, we proceed as follows.

For each face $f_\partial$ in $\mc{F}_\partial$, we let $\{T_1, T_2, N\}$ be 
an oriented orthonormal frame, with $N$ the exterior normal to the said 
face in $K$. The barycentric subdivision of $f_\partial$ will have 
twelve edges $e^j_{f_\partial}$, $j=1, \ldots, 12$
and six faces $f^k_{f_\partial}$, $k=1, \ldots, 6$. For notational
convenience, we denote the linear combination 
$\sum_{j=1}^{12} c_{e^j_{f_\partial }} 
W_{e^j_{f_\partial}} + \sum_{k=1}^6 c_{f^k_{f_\partial }} W_{f^k_{f_\partial}}$
by $W_{c_{e^j_{f_\partial }}e^j_{f_\partial}} + 
W_{c_{f^k_{f_\partial }}f^k_{f_\partial}}$.

Since for any edge $e$ the matrix of component derivatives $\nabla W_{e^j}$ is 
antisymmetric, by (\ref{bs}) and the symmetries of the tensor of 
elastic constants $W$, we have that 
$(\sigma(\nabla W_{e^j_{f_\partial}}) N)^\alpha =  
W^{\alpha \hspace{1mm} r}_{\mbox{}\hspace{1.5mm}i \hspace{1.5mm}\beta}
\, \partial_r W_{c_{e^j_{f_\partial}}e^j_{f_\partial}}^\beta N^i =0$. 
Hence, we discretize the local condition in (\ref{nblg}) by requiring that 
\begin{equation} \label{e27}
\begin{array}{rcl}
\< \partial _{T_1}(W_{c_{e^j_{f_\partial }}e^j_{f_\partial}} +
W_{c_{f^k_{f_\partial }}f^k_{f_\partial}}),W'(\BOne)N\>+
\< T_1, 
\sigma(\nabla W_{c_{f^k_{f_\partial}} f^k_{f_\partial}}) N) \>
& =  & 0 \, , \\
\< \partial _{T_2}(W_{c_{e^j_{f_\partial }}e^j_{f_\partial}} +
W_{c_{f^k_{f_\partial }}f^k_{f_\partial}}),W'(\BOne)N\>+
\< T_2, 
\sigma(\nabla W_{c_{f^k_{f_\partial}} f^k_{f_\partial}}) N) \>
& =  & 0 \, , 
\end{array}
\end{equation}
where $(W'(\BOne)N)^{\alpha}  = {\displaystyle
W^{\alpha \hspace{0.5mm} j}_{\mbox{}\hspace{1.4mm}i \hspace{1mm}j}N^i}$,
and the pairings are in the sense of $L^2$ of the boundary face 
$f_\partial$. By the coercivity of the tensor $W$, at most one, if at all,
of the last summands on the left in these expressions can be zero.

This $2\times 18$ system of homogeneous equations couples the indicated 
barycentric boundary coordinates of $U$ associated to the face $f_\partial \in
\mc{F}_\partial$. In addition to that,  
the blocks associated to adjacent faces in $\mc{F}_\partial$ couple between 
them the two barycentric boundary edge coordinates of $U$ associated to the one
common edge of these faces. 

Since $|\mc{E}'_\partial| =2|\mc{E}_\partial|+6| \mc{F}_\partial|$ 
and $| \mc{F}'_\partial|=6| \mc{F}_\partial|$, the procedure carried 
out over all the faces of $\mc{F}_\partial$ produces 
an underdetermined system of $2|\mc{F}_\partial |$ homogeneous equations in the 
$(2| \mc{E}_\partial| + 6| \mc{F}_\partial|) + 6| \mc{F}_\partial|$ 
boundary coordinate unknowns $\left( \begin{array}{c}
        c_{e_{f_\partial}} \\ c_{f_{f_\partial}}
        \end{array}
\right)
$
of $U$. 

By the symmetries and coercivity of $W$, the
rank of the $2\times 18$ block associated to each  
$f_{\partial }\in \mc{F}_\partial$ is either two, or one, 
generically the former, and equals the
rank of the $2 \times 6$ subblock in it that
involves the barycentric faces only.
Thus, the equations in this block 
can be used to express two, or one, of the coordinates 
$c_{f^k_{ f_\partial}}$ as a linear combination of the remaining edge 
and face barycentric boundary coordinates involve in it, unaffected
by the extra coupling of adjacent face equations mentioned above.
We carry out the row reduction of the block, and that of the entire system,
accordingly.

By a suitable reordering of the basis elements, we may write 
the row reduced matrix of the entire system as 
$(-C \; \; \BOne\, )$, where $C$ is a block whose number of rows and columns 
are bounded above by $2| \mc{F}_\partial |$ and $2|\mc{E}_{\partial}|+
11| \mc{F}_\partial|$, 
respectively. We decompose the boundary faces in $\mc{F}'$ accordingly, 
$\mc{F}'_\partial = \mc{F}'_{\partial^I}
+\mc{F}'_{\partial^B}$, so that 
$c_{e_\circ}W_{e_\circ}+ c_{e_\partial}W_{e_\partial}+
 c_{f_\circ}W_{f_\circ}+ c_{f_{\partial^I}}W_{f_{\partial^I}}
 +c_{f_{\partial^B}}W_{f_{\partial^B}}
\in L_1^2(K';\mb{R}^3)$ satisfies (\ref{e27}) if, and only if,
$c_{f_{\partial^B}}=C \left( \begin{array}{c}
        c_{e_{\partial }} \\ 
        c_{f_{\partial^I}} 
        \end{array} \right)$, 
and define $L^2_{1,b}(K';\mb{R}^3)$ as such a subspace:
\begin{equation} \label{e30bis} 
L^2_{1,b}(K')  = \{ 
c_{e_\circ}W_{e_\circ}+ c_{e_\partial}W_{e_\partial}+
 c_{f_\circ}W_{f_\circ}+ c_{f_{\partial^I}}W_{f_{\partial^I}}
 +c_{f_{\partial^B}}W_{f_{\partial^B}}: 
\; c_{f_{\partial^B}}=C \left( \begin{array}{c}
        c_{e_{\partial }} \\
        c_{f_{\partial^I}} 
        \end{array} \right)
\} \, . 
\end{equation}
If the row reduced matrix of the system of boundary conditions
(\ref{e27}) were not to have any null rows for the $2\times 6$ subblocks, 
we would have that  
${\rm dim}\, L^2_{1,b}(K';\mb{R}^3)=|\mc{E}'|+|\mc{F}'|-2| \mc{F}_\partial |=
|\mc{E}'| +|\mc{F}'_{\circ}|+4| \mc{F}_{\partial}|$. 
Otherwise ${\rm dim}\, L^2_{1,b}(K';\mb{R}^3)=|\mc{E}'|+
|\mc{F}'|-\left( 2(| \mc{F}_\partial |-r_n)+
r_n\right)=|\mc{E}'|+
 |\mc{F}'_{\circ}|+4| \mc{F}_{\partial}|+r_n$, where $r_n$ is the
number of null rows of the $2 \times 6$ subblocks. If 
$C=\left( c_{f_{\partial^B},e_{\partial}} \; c_{f_{\partial^B},
f_{\partial^I}} \right)$, 
the set $\mc{B}_{L^2_{1,b}}=\{ W_s \}_{W_s \in L^2_{1,b}}:=
\{ W_{e_{\circ}}, W_{f_\circ}, W_{e_{\partial}} + 
c_{f_{\partial^B},e_{\partial}}W_{f_{\partial^B}},   
W_{f_{\partial^I}}+c_{f_{\partial^B},f_{\partial^I}} 
W_{f_{\partial^B}}\}_{e_\circ \in \mc{E}'_\circ, f_\circ \in \mc{F}'_\circ,
e_{\partial}\in \mc{E}'_{\partial}, f_{\partial^I}\in \mc{F}'_{
\partial^I}}$ is a basis for $L^2_{1,b}(K';\mb{R}^3)$. 

\subsection{The discretized equation} \label{tde}
The dualization of the term 
$A_{ij}^{\alpha \beta}(\BOne)\partial_i \partial_j u=
W^{i\hspace{1.5mm}j}_{\phantom{i}\alpha \hspace{1.5mm} \beta}
\partial_i \partial_ju$ is accomplished by writing this operator in divergence
form. We obtain 
$$
\< A_{ij}^{\alpha \beta}(\BOne)\partial_i \partial_j u,u\>=
-\< A_{ij}^{\cdot \, \beta}(\BOne) \partial_j u,\partial_i u\>+
\< \sigma(\nabla u) N, u\mid_{\partial \Omega}\>_{L^2(\partial \Omega)} \, ,  
$$
where the stress tensor boundary term follows by using
(the middle relation in) (\ref{bs}). It is then clear how this term is
discretized. Notice that by the antisymmetry of
$\nabla W_e$, the induced bilinear form that these two terms associate
with pairs of edge coefficients of any type is zero, and by the symmetries
of the tensor $W$, the same of true for any edge-face pair also.   

If $h=h(u)$ solves the boundary value problem (\ref{hb}), we have that  
$$
\< \nabla h, u\> = - (\< \sigma(\nabla u)N,N\>+
\<W'(\BOne)N,N\>\partial_N u, N\>)\<u, N\>_{L^2(\partial \Omega)} - 
\< h, {\rm div}\, u\> \, .
$$
The discretization 
of the boundary term on the right is straightforward
for $U\in L^2_1(K';\mb{R}^3)$. For the discretization of the last pairing,
we observe that
$-\Delta h= L_1(D)\partial_1 u^1 + L_2(D)\partial_2 u^2 + L_3(D)\partial_3u^3$,
where   
\begin{equation}
\begin{array}{rcl}
L_1(D) & = &  W^{1\hspace{1.5mm}1}_{\phantom{1}1\hspace{1.5mm} 1}\partial_1^2
+(2W^{2\hspace{1.5mm}2}_{\phantom{2}1\hspace{1.5mm}1}+
W^{2\hspace{1.5mm}1}_{\phantom{2}2\hspace{1.5mm}1})\partial_2^2+
(2W^{3\hspace{1.5mm}3}_{\phantom{3}1\hspace{1.5mm}1}
+W^{3\hspace{1.5mm}1}_{\phantom{3}3\hspace{1.5mm}1})\partial_3^2\,  ,
 \vspace{1mm} \\
L_2(D) & = &  
(2W^{1\hspace{1.5mm}1}_{\phantom{2}2\hspace{1.5mm}2}+
W^{1\hspace{1.5mm}2}_{\phantom{2}1\hspace{1.5mm}2})\partial_1^2+
W^{2\hspace{1.5mm}2}_{\phantom{2}2\hspace{1.5mm} 2}\partial_2^2 +
(2W^{2\hspace{1.5mm}2}_{\phantom{2}3\hspace{1.5mm}3}
+W^{3\hspace{1.5mm}2}_{\phantom{3}3\hspace{1.5mm}2})\partial_3^2 \, , 
\vspace{1mm}\\
L_3(D) & = &  
(2W^{1\hspace{1.5mm}1}_{\phantom{1}3\hspace{1.5mm}3}+
W^{1\hspace{1.5mm}3}_{\phantom{1}1\hspace{1.5mm}3})\partial_1^2+
(2W^{2\hspace{1.5mm}2}_{\phantom{2}3\hspace{1.5mm}3}+
+W^{2\hspace{1.5mm}3}_{\phantom{3}2\hspace{1.5mm}3})\partial_2^2
+W^{3\hspace{1.5mm}3}_{\phantom{3}3\hspace{1.5mm}3}\partial_3^2\, .
\end{array}
\end{equation}
Each of these elliptic nonpositive operators is weakly equivalent to the
Laplacian. If $l_{L_i(D)}$ is the average of the coefficients of 
$L_i(D)$, for numerical purposes, we choose to approximate 
$L_i(D)$ by $l_{L_i(D)}\Delta$, so we will have
$-\Delta h \cong l_{L_1(D)}\Delta \partial_1 u^1 + l_{L_2(D)}\Delta
\partial_2 u^2 + l_{L_3(D)}\Delta \partial_3u^3$. This induces 
the $L^2(\Omega)$ approximation to $h$ given by the weighted divergence   
\begin{equation} \label{wedi}
-h(u) \cong l_{L_1(D)}\partial_1 u^1 + l_{L_2(D)}
\partial_2 u^2 + l_{L_3(D)}\partial_3u^3\, , 
\end{equation}
which we use, and discretize it to $h(U)$. Notice that
for neoHookian materials, $L_1(D)=L_2(D)=L_3(D)=\Delta$, in which case we
have that $-h(u)={\rm div}\, u$ exactly.

As for the last term, by the divergence condition in (\ref{nblg}), 
the dualization of
$\lambda \nabla {\rm div}\, u$ have zero boundary contribution.
We have that
$$
\lambda \< \nabla {\rm div} \, u, u\> = -\lambda \< {\rm div}\, u,
{\rm div}\, u\> \, , 
$$
which can be discretized in the obvious manner. 
Notice that ${\rm div}\, W_e  =0$ in $L^2(K')$, so the induced bilinear 
form that this term associates with edge coefficients of any type, boundary
or interior, is zero.
Since the parameter $\lambda$ is introduced in (\ref{new2})
just to feel the effect of linearizing $J(\eta(t))$ when $J(\eta(t))$
does not necessarily satisfy (\ref{eq1}), it is reasonable to take
as $\lambda$ the value 
\begin{equation} \label{feel}
\lambda:=\lambda_L=1  \, .
\end{equation}
We show the naturality of this choice below by carrying a simulation 
with $\lambda$ equal to the average of the $l_{L_i(D)}$s also, which for the
$l_{L_i}$s of our simulations, is a substantially larger value, and which makes 
this term of the same magnitude as the previous two, possibly leading to 
undesirable cancellations.

If over dots stand for time derivatives, the coefficients of $U$ in
(\ref{asb}) are solutions of the equation 
\begin{equation} \label{e24}
\rho I_{\circ, \partial}
 \left( \begin{array}{c}
\ddot{c}_{e'_\circ} \\ \ddot{c}_{e'_{\partial}} \\ 
\ddot{c}_{f'_\circ} \\ \ddot{c}_{f'_{\partial}} 
\end{array} 
\right)= K_{\circ , \partial }
 \left( \begin{array}{c}
c_{e'_\circ} \\ c_{e'_{\partial}}\\ 
c_{f'_\circ} \\ c_{f'_{\partial}} 
\end{array} 
\right) +
\left( \begin{array}{c}
\< F, W_{e_\circ}\>  \\ \< F, W_{e_{\partial}} \> \\
\< F, W_{f_\circ}\>  \\ \< F, W_{f_{\partial}} \> 
\end{array} 
 \right) \, ,
\end{equation}
where the block decomposition of the matrix $I_{\circ , \partial }$ of inner 
products of the basis elements is explicitly given by 
$$
I_{\circ , \partial}
= \left( \begin{array}{cccc}
\<W_{e_\circ},W_{e'_\circ}\> &  \< W_{e_\circ},W_{e'_{\partial}}\> & 
\<W_{e_\circ},W_{f'_\circ}\> &  \< W_{e_\circ},W_{f'_{\partial}}\> \\
\< W_{e_{\partial}},W_{e'_\circ}\> & \<W_{e_{\partial}},W_{e'_{\partial} }\> &
\< W_{e_{\partial}},W_{f'_\circ}\> & \<W_{e_{\partial}},W_{f'_{\partial} }\> \\
\<W_{f_\circ},W_{e'_\circ}\> &  \< W_{f_\circ},W_{e'_{\partial}}\> & 
\<W_{f_\circ},W_{f'_\circ}\> &  \< W_{f_\circ},W_{f'_{\partial}}\> \\
\< W_{f_{\partial}},W_{e'_\circ}\> & \<W_{f_{\partial}},W_{e'_{\partial} }\> &
\< W_{f_{\partial}},W_{f'_\circ}\> & \<W_{f_{\partial}},W_{f'_{\partial} }\> 
\end{array}
\right) \, ,
$$
and the block decomposition of the matrix $K_{\circ , \partial }$ is of
the form
$$
K_{\circ , \partial}
= \left(  \begin{array}{cccc}
\mb{O}_{e_\circ,e'_\circ} &  \mb{O}_{e_\circ,e'_{\partial}} &
\mb{O}_{e_\circ,f'_\circ} &  \mb{O}_{e_\circ,f'_{\partial}} \\
\mb{O}_{e_{\partial},e'_\circ} & 
\mb{O}_{e_{\partial},e'_{\partial} } &
\mb{O}_{e_{\partial},f'_\circ} & K_{e_{\partial},f'_{\partial} } \\
\mb{O}_{f_\circ,e'_\circ} &  \mb{O}_{f_\circ,e'_{\partial}} &
K_{f_\circ,f'_\circ} &  K_{f_\circ,f'_{\partial}} \\
\mb{O}_{f_{\partial},e'_\circ} & K_{f_{\partial},e'_{\partial} } &
K_{f_{\partial},f'_\circ} & K_{f_{\partial},f'_{\partial} }
\end{array}
\right) \, , 
$$
where $\mb{O}_{ij}$ stands for the $ij$th entry of the zero matrix, 
the remaining blocks associated with boundary edge elements given by  
$$
\begin{array}{rcl}
K_{e_\partial ,f'_{\partial}} & = & 
\frac{1}{2}  \<\sigma(\nabla W_{f'_\partial})N, W_{e_\partial}\>
-\frac{1}{2}  \<\sigma(\nabla W_{f'_\partial})N,N\> \< W_{e_{\partial}},N\>
\\ & &  
-\frac{1}{2}
\< W'(\BOne) N,N\>\left( 
\<\partial_N W_{f_{\partial}'}, N\>\< W_{e_{\partial}},N\>
+\<\partial_N W_{e_{\partial}}, N\>\< W_{f_{\partial}'},N\>\right) \, , \\
K_{f_\partial ,e'_{\partial}} & = & 
\frac{1}{2}  \<\sigma(\nabla W_{f_\partial})N, W_{e'_\partial}\>
-\frac{1}{2}  \<\sigma(\nabla W_{f_\partial})N,N\> 
\< W_{e'_{\partial}},N\> \vspace{1mm} \\ & & 
-\frac{1}{2}
\< W'(\BOne) N,N\>(\<\partial_N W_{e_{\partial}'}, N\>\< W_{f_{\partial}},N\>
+\<\partial_N W_{f_{\partial}}, N\>\< W_{e_{\partial}'},N\>) \, , 
\end{array}
$$
the pairings on the right here in the sense of $L^2(\partial \Omega)$, and  
$$
\begin{array}{rcl}
K_{f_\circ,f'_\circ} & = & 
-\< \partial_i W^{\alpha}_{f_\circ}, A^{\alpha \beta}_{ij}\partial_j
W^\beta_{f'_\circ}\>
-\lambda_L \< \partial_i W^i_{f_\circ}, \partial_j W^j_{f'_\circ}\> +\\ &  & 
\hspace{15mm}
\frac{1}{2}(\< l_{L_i(D)}\partial_i W^i_{f_\circ}, \partial_j W^j_{f'_\circ}\>
+\< l_{L_i(D)}\partial_i W^i_{f'_\circ}, \partial_j W^j_{f_\circ}\>)\, , 
\vspace{1mm} \\
 K_{f_\circ,f'_{\partial}} & = & 
-\< \partial_i W^{\alpha}_{f_\circ}, A^{\alpha \beta}_{ij}\partial_j
W^\beta_{f'_{\partial}} \>
-\lambda_L \< \partial_i W^i_{f_\circ}, \partial_j W^j_{f'_\partial}\> +\\ & &
\hspace{15mm}
\frac{1}{2}(\< l_{L_i(D)}\partial_i W^i_{f_\circ}, \partial_j W^j_{f'_\partial
}\>
+\< l_{L_i(D)}\partial_i W^i_{f'_\partial}, \partial_j W^j_{f_\circ}\>)\, ,
\vspace{1mm} \\
K_{f_{\partial},f'_\circ} & = & 
-\< \partial_i W^\alpha_{f_{\partial}},
A^{\alpha \beta}_{ij}\partial_j W^\beta_{f'_\circ}\>
-\lambda_L \< \partial_i W^i_{f_\partial}, \partial_j W^j_{f'_\circ}\> +\\ & & 
\hspace{15mm}
\frac{1}{2}(\< l_{L_i(D)}\partial_i W^i_{f_\partial}, \partial_j W^j_{f'_\circ
}\>
+\< l_{L_i(D)}\partial_i W^i_{f'_\circ}, \partial_j W^j_{f_\partial}\>)\, ,
\vspace{1mm} \\ 
 K_{f_{\partial},f'_{\partial} } & = & 
 -\< \partial_i W^\alpha_{f_{\partial}},
A^{\alpha \beta}_{ij} \partial_j W^\beta_{f'_{\partial} }\>  
-\lambda_L \< \partial_i W^i_{f_\partial }, \partial_j W^j_{f'_\partial}\>\\ & &
+\frac{1}{2}(\< l_{L_i(D)}\partial_i W^i_{f_\partial}, \partial_j 
W^j_{f'_\partial }\>
+\< l_{L_i(D)}\partial_i W^i_{f'_\partial}, \partial_j W^j_{f_\partial}\>)
\vspace{1mm} + 
B(W_{f_{\partial}}, W_{f_{\partial}'})\,  , 
\end{array}
$$
where the boundary term $B(W_{f_{\partial}}, W_{f_{\partial}'})$
is given by the sum of $L^2(\partial \Omega)$ pairings 
$$
B(W_{f_{\partial}}, W_{f_{\partial}'}) =
\< \sigma(\nabla W_{f_{\partial}'})N,W_{f_{\partial}}\> -
(\< \sigma(\nabla W_{f_{\partial}'})N,N\> +\< W'(\BOne) N,N\>\<\partial_N
W_{f_{\partial}'}, N\>)\< W_{f_{\partial}},N\> \, .
$$

The matrices in (\ref{e24}) are sparse, but their sparsity is lower
than the sparsity of the square block associated to faces only. (This latter
number agrees with the sparsity of the matrices in \cite[system (31)]{sim},
which are, up to the definition of the $K$ face blocks here, the same.) 
Indeed, for the entry of an edge block that a pair of edges $e,e'$ defines to 
be zero, it is sufficient that ${\rm St} e\cap {\rm St} e'= \emptyset$, while
for the entry of a cross block that a pair $e, f'$ of edge and face
defines to be zero, it is sufficient that
${\rm St} e\cap {\rm St} f'= \emptyset$. But there are fewer pairs $e,e'$,
 or $e,f'$, satisfying these conditions  
relative to their total than pair of faces $f,f'$ relative to their total
satisfying the condition ${\rm St} f\cap {\rm St} f'= \emptyset$, which 
suffices for the entries in the matrices defined by $f,f'$ to vanish. 

The system (\ref{e24}) is symmetric. We notice that the relatively 
insignificant presence
of edge terms in $K_{\circ, \partial}$ reflects the fact that edge elements
$W_e$  are used to approximate the gradient component of $U$, 
a vector field that the equation of motion tries to bring closer to 
a divergence-free vector field starting from one that is already 
relatively close. Thus, edges play a lesser role in finding $U$
than that played by faces, and edges enter into the definition of
$K_{\circ,\partial}$ only when interacting with a face.

With the vector of coefficients given by the solution to (\ref{e24}), 
(\ref{asb}) produces a numerical approximation to our solution of (\ref{vp2}).
The eigenvalues of $-(\rho I_{\circ, \partial})^{-1} K_{\circ,
\partial}$ and the corresponding frequencies they induce approximate  
natural vibration frequencies of the body. Positive eigenvalues lead to 
vibrations that decay exponentially fast in time. Negative eigenvalues lead
to undamped vibration modes. By resonance, any one of the latter
induces oscillatory motions within the body when this is subjected to an 
external sinusoidal pressure wave force $F=F_{f}$ of frequency $f$ close 
to the frequency of the wave mode.

Since the eigenvectors of $-(\rho I_{\circ, \partial})^{-1} 
K_{\circ, \partial}$ do not necessarily satisfy the boundary condition
(\ref{e27}), we take the waves they induce 
as {\it coarse} approximations to the vibration patterns of the 
body, and call {\it coarse resonance waves} the waves produced by resonance 
for frequencies close to the frequencies they intrinsically have.

By construction, the global divergence condition in (\ref{nblg}) 
is satisfied by any of the waves above. We derive the {\it fine} 
approximations to the solution of (\ref{vp2}) by incorporating into these
waves the remaining local conditions, which in its discrete form, are 
given by the system (\ref{e27}). We discretize the solution $U$ now over 
the space $L^2_{1,b}(K';\mb{R}^3)\subset L^2_1(K';\mb{R}^3)$ in (\ref{e30bis}) 
instead, using the natural basis $\mc{B}_{L^2_{1,b}}$ for this space that 
the splitting $\mc{F}'_{\partial}= \mc{F}'_{\partial^I}+ \mc{F}'_{\partial^B}$
leads to, and expressing the discretized solution into blocks accordingly, 
\begin{equation} \label{asbf}
U = c_{e_\circ}W_{e_{\circ}} + 
c_{e_{\partial}^B}( W_{e_\partial}+
c_{f_{\partial^B},e_{\partial}} W_{f_{\partial^B}})
+ c_{f_\circ} W_{f_\circ}+
c_{f_{\partial^I}^B}(W_{f_{\partial^I}}+
c_{f_{\partial^B},f_{\partial^I}}W_{f_{\partial^B}}) \, , 
\end{equation}
where $C=(c_{f_{\partial^B},e_{\partial}} \; 
c_{f_{\partial^B},f_{\partial^I}}):=
(C_{e_\partial} \; C_{f_{\partial^I}})$ 
is the $C$ block of the row reduced 
matrix of conditions (\ref{e27}). We obtain a square system    
of differential equations for the coefficients 
\begin{equation} \label{e24f}
\rho I_{\circ, \partial^B}
 \left( \begin{array}{c}
\ddot{c}_{e_\circ} \\ \ddot{c}_{e_{\partial}^B} \\
\ddot{c}_{f_\circ} \\ \ddot{c}_{f_{\partial^I}^B}
\end{array}
\right)= K_{\circ , \partial^B }
 \left( \begin{array}{c}
c_{e_\circ} \\ c_{e_{\partial}^B}\\
c_{f_\circ} \\ c_{f_{\partial^I}^B}
\end{array}
\right) +
\left( \begin{array}{c}
\< F, W_{e_\circ}\>  \\ \< F, W_{e_{\partial}} +
c_{f_{\partial^B},e_{\partial}}W_{f_{\partial^B}} \> \\
\< F, W_{f_\circ}\>  \\ 
\< F, W_{f_{\partial^I}}+ c_{f_{\partial^B},f_{\partial^I}}
W_{f_{\partial^B}}  \>
\end{array}
 \right) \, ,
\end{equation}
which is, of course, very closely related to (\ref{e24}) but now incorporates
the additional splitting induced by the decomposition 
$\mc{F}'_\partial= \mc{F}'_{\partial^I}+\mc{F}'_{\partial^B}$ of the 
boundary faces. 
For instance, the blocks $I^f_{\circ, \partial^B}$ and 
$K^f_{\circ, \partial^B}$ in $I_{\circ, \partial^B}$ and 
$K_{\circ, \partial^B}$ that are associated to pair of faces only are given by 
$$
\begin{array}{rcl}
I^f_{\circ,\partial^B} & = & \left( \! \! \begin{array}{cc}
I_{{f_\circ},{f'_\circ}} & I_{{f_\circ},{f'_{\partial^I}}}
+I_{{f_\circ},{f_{\partial^B}'}} C^t_{f_{\partial^I}} \\ 
I_{f_{\partial^I},{f'_{\circ }}}+
C_{f_{\partial^I}}I_{f_{\partial^B},{f'_{\circ}}} & 
C_{f_{\partial^I}} I_{{f_{\partial^B}},{f'_{\partial^B}}}
C^t_{f_{\partial^I}}+
C_{f_{\partial^I}} I_{{f_{\partial^B}},{f'_{\partial^I}}}
+
I_{{f_{\partial^I}},{f'_{\partial^B}}} C^t_{f_{\partial^I}} +
I_{{f_{\partial^I}}, {f_{\partial^I}}'}
\end{array}
\! \! \right) \, , \vspace{1mm} \\
K^f_{\circ, \partial^B} & = & \left( \begin{array}{cc}
K_{f_\circ,f'_\circ} &  
K_{f_\circ,f'_{\partial^I}}+ K_{f_\circ,{f_{\partial^B}'}} 
C^t_{f_{\partial^I}} \\ 
K_{f_{\partial^I},f'_\circ} +C_{f_{\partial^I}}K_{f_{\partial^B},f'_\circ} & 
C_{f_{\partial^I}} K_{f_{\partial^B},f'_{\partial^B} }C^t_{f_{\partial^I}}+
C_{f_{\partial^I}} K_{{f_{\partial^B}},{f_{\partial^I}'}}
+
K_{{f_{\partial^I}},{f'_{\partial^B}}}
C^t_{f_{\partial^I}} +
K_{f_{\partial^I},{f'_{\partial^I}}} 
\end{array}
\right) \, , 
\end{array}
$$
respectively, where the $I,K$ entries on the right are defined by the 
expressions given in (\ref{e24}). (The remaining blocks in (\ref{e24f}) 
have a similar description, incorporating the role that 
$C_{e_\partial}$ plays.) The system so produced is symmetric.

Notice that the local boundary conditions (\ref{e27}) bring about relations
into the system through the boundary edges; the blocks in 
$K_{\circ,\partial^B}$  associated to pairs of boundary edges, or a boundary 
edge and an interior face, are now nontrivial in comparison with the 
analogous entries for the matrix $K_{\circ, \partial }$ in the
system (\ref{e24}). Notice also that the bottom right blocks for the 
matrices $I^f_{\circ, \partial^B}$, and $K^f_{\circ,\partial^B}$, respectively,
are no longer diagonal, as was the case of the corresponding blocks in 
the matrices of (\ref{e24}). 

With the coefficients given by the solution to (\ref{e24f}), the vector
field (\ref{asbf}) yields an approximation to the solution of
(\ref{vp2}) that satisfies the boundary conditions (\ref{nblg}).  
The {\it fine} approximations to the vibration patterns of the body are  
those waves induced by the eigenvectors of the matrix  
$-(\rho I_{\circ, \partial^B})^{-1} K_{\circ, \partial^B}$. They
satisfy the local boundary conditions in (\ref{nblg}) by construction, and 
therefore, correspond to waves producing curves in $L^2_{1,b}(K';\mb{R}^3)$.
The {\it fine resonance waves} are those associated to the negative
eigenvalues of this matrix, generated by resonance 
when $\Omega$ is acted on by the external sinusoidal pressure wave force 
$F=F_{f}$ of frequency $f$ close to the intrinsic frequency of the fine 
wave modes.

\subsection{The full fledged first iterate $(\zeta, \dot{\zeta})=
\mc{M}_F(\BOne,w)$} 
If the initial conditions $(u_0,u_1)$ of $u(t)$ in (\ref{vp2}) are 
compatible with the initial conditions $(\BOne , w)$ for (\ref{new2}), 
our algorithm would produce a numerical solution $U=U(t)$ of (\ref{e24f}) 
with initial condition $(U_0,U_1)$ that is the discretization 
of $(u_0=w, u_1)$, and the pair $(U(t), V(t))$, $V(t)=\dot{U}(t)$, would be 
the numerical version of $S_F(\BOne, w)$ discretized over the spaces that 
we introduced for the purpose. This pair, together with 
the discretized pressure force $F$ in the right hand side 
of (\ref{e24f}), could then be used to produce a discretization over these 
spaces of the right hand side of the system  (\ref{imp}), a suitable
numerical solution of which would then be a numerical version of the 
full fledged first iterate 
$(\zeta , \dot{\zeta})= \mc{M}_F(\BOne, w)$ of the map (\ref{nohm}). 

The component $\zeta$ of one such is a $C^1$ path of embeddings that, 
for each $t$, is a linear combinations of all the modes of vibrations of 
$U(t)$, damping and oscillatory. We present an algorithm to compute 
it, the algorithm being a function of $U(t)$, and the discretized $F$. We do 
not make use of it in our work, but we can conceive of situations 
where, in spite of the technical difficulties given the complexity,
we could use this numerical $\zeta(t)$ to analyze the 
vibrations of the body in further detail from the details provided here 
(in particular, if we were to
concentrate our attention on the vibrations at frequencies nearby a 
predetermined one), and so we think it useful to have a way of computing it
readily available, in case these details are of importance to obtain.

We recall that 
$$
\gamma^2 \rho \zeta - A_{\lambda}(\zeta, \dot{\zeta})
=-\rho v(t)+ \int_0^t F(s) ds 
+\gamma^2 \rho \left( \BOne + \int_0^t u(s) ds \right)\, ,  
$$
where $A_\lambda(\zeta, \dot{\zeta})=A_{\lambda=1}(\zeta,\dot{\zeta})$ is 
given by
$$
A_1(\zeta, \dot{\zeta})
=A_{ij}^{\alpha \beta }(D\zeta )\partial _{i}
\partial _{j}\zeta ^{\beta } +
J(\zeta(t))\nabla_\zeta q + J(\zeta(t)) \nabla_\zeta
J(\zeta(t) ) \, ,
$$
$q$ being the solution to the boundary value problem
$$ 
\begin{array}{rcl}
{\rm div}_{\zeta} J(\zeta) \nabla_{\zeta }q
& = & -{\rm div}_{\zeta }A_{ij}^{\alpha \beta }(D\zeta )\partial _{i}
\partial _{j}\zeta ^{\beta } +\rho \,
{\rm trace}(D_{\zeta } \dot{\zeta })^{2}\, ,  \vspace{1mm} \\
q\mid _{\partial \Omega }\; & = &  \;-{\displaystyle \frac{\langle W'(D\zeta) N
 , \nu \circ \zeta  \rangle }{J^{b}(\zeta )}} \; ,
\end{array}
$$
and that $\zeta$ satisfies the conditions
$$
\text{$\< W'(D\zeta ) N , \partial_T \zeta\>=0$, \quad $J(\zeta)=1$, 
\quad
on $\partial \Omega $.} 
$$ 

We freeze $\zeta(t)=\BOne$ in the operators $J(\zeta(t))\nabla_{\zeta}$, 
$D_{\zeta}$, and ${\rm div}_{\zeta}$ in these expressions, 
as well as in the terms $J^b(\zeta)$ and $J({\zeta})$ in the boundary 
conditions for $q$, and $\zeta$, respectively. Accordingly, we replace 
$D\zeta$ by $D\zeta =\BOne$ in $A_{ij}^{\alpha \beta }(D\zeta )$ and 
$W'(D\zeta)$, and set $\dot{\zeta}=u(t)$ in the quadratic 
trace term in the right of the interior equation for $q$. The term 
$J(\zeta(t)) \nabla_\zeta J (\zeta(t))$ in the right side of 
the defining expression for $A_1 (\zeta, \dot{\zeta})$ becomes 
$\nabla J(\zeta(t))$, and since $J(\zeta(t))$ is an ordered cubic in the 
components of $\nabla \zeta$, we freeze further 
$\zeta(t)=\BOne$ in the first two factors of this cubic, thus transforming
 $\nabla J(\zeta(t))$ into $\nabla {\rm div} \, \zeta$. 
By (\ref{bs}), all of this results into the equation
\begin{equation} \label{fif}
\gamma^2 \rho  \zeta - \mc{A}(\zeta) =-\rho v(t)+ \int_0^t F(s) ds
+\gamma^2 \rho \left( \BOne + \int_0^t u(s) ds \right) 
\end{equation}
for $\zeta$, where $\mc{A}$ is the operator
\begin{equation} \label{ope}
\mc{A}(\zeta)
=W^{i\hspace{1.5mm}j}_{\phantom{i}\alpha \hspace{1.5mm} \beta} \partial _{i}
\partial _{j}\zeta ^{\beta } + \nabla q + \nabla {\rm div} \, \zeta \, ,
\end{equation}
$q$ is now the solution of the boundary value problem
\begin{equation} \label{linq}
\begin{array}{rcl}
{\rm div}\, \nabla q= \Delta q
& = & -{\rm div}W^{i\hspace{1.5mm}j}_{\phantom{i}\alpha \hspace{1.5mm} \beta} 
\partial _{i}
\partial _{j}\zeta ^{\beta } +\rho \,
{\rm trace}(\nabla u)^{2}\, ,  \vspace{1mm} \\
q\mid _{\partial \Omega }\; & = &  \;-{\displaystyle 
\langle W^{i\hspace{1.5mm}j}_{\phantom{i}\alpha \hspace{1.5mm} j} N^\alpha
 , \nu \circ \zeta  \rangle } \; ,  
\end{array}
\end{equation}
and $\zeta$ satisfies the conditions 
\begin{equation} \label{fibc}
\text{$\< W^{i\hspace{1.5mm}j}_{\phantom{i}\alpha \hspace{1.5mm} j}N^\alpha,
\partial_T \zeta\>=0$  on $\partial \Omega$.} 
\end{equation}

By the coercivity of the stored energy function, if $\gamma$ is chosen to be
sufficiently large, the operator 
$$
\gamma^2 \rho \, \BOne - \mc{A}: H^2(\Omega; \mb{R}^3) 
\rightarrow L^2(\Omega; \mb{R}^3) 
$$
is invertible, and (\ref{fif}) can be solved for $\zeta$, with $\zeta
\in H^2(\Omega;\mb{R}^3)$. We find this $\zeta$ by solving
numerically a weak version of this equation, with the boundary conditions
(\ref{fibc}) enforced upon the solution. 

We discretize the sought after solution $\zeta$ over the space 
$L^2_1(K';\mb{R}^3)$, 
\begin{equation} \label{dfif}
Z= d_{e}W_{e}+ d_{f}W_{f}=d_{e_\circ}W_{e_\circ}+ d_{e_\partial}W_{e_\partial}+
 d_{f_\circ}W_{f_\circ}+ d_{f_\partial}W_{f_\partial}\, . 
\end{equation}
Proceeding as in \S \ref{sbc}, we make $Z$ satisfy the conditions  
(\ref{fibc}) by requiring that 
\begin{equation}\label{dfibc} 
\begin{array}{rcl}
\< \partial _{T_1}(W_{d_{e^j_{f_\partial }}e^j_{f_\partial}} +
W_{d_{f^k_{f_\partial }}f^k_{f_\partial}}),W'(\BOne)N\> & =  & 0 \, , \\
\< \partial _{T_2}(W_{d_{e^j_{f_\partial }}e^j_{f_\partial}} +
W_{d_{f^k_{f_\partial }}f^k_{f_\partial}}),W'(\BOne)N\> & =  & 0 \, ,
\end{array}
\end{equation}
for each face $f_\partial$ in $\mc{F}_\partial$. Here,
$e^j_{f_\partial}$, $j=1, \ldots, 12$ and 
$f^k_{f_\partial}$, $k=1, \ldots, 6$, are the twelve edges and six faces in
the barycentric subdivision of the face $f_\partial$, respectively, 
with associated coefficients $d_{e^j_{f_\partial }}$, and 
$d_{f^k_{f_\partial }}$, 
and $\{T_1, T_2, N\}$ is an oriented orthonormal frame of the tangent space,  
with $N$ the exterior normal to the said
face in $K$. 

The whole vector of coordinates $\left( \begin{array}{c}
        d_{e} \\ d_{f}
        \end{array}
\right)
$
is found as the solution of the linear system that results from the weak 
formulation
of the elliptic equation (\ref{fif}), subject to the constrains (\ref{dfibc}).

The dualizations of the first and last terms on the right of (\ref{fif}) are 
straightforward, an almost verbatim repetition of the analogous arguments
in \S \ref{tde}:
$$
\begin{array}{rcl} 
\< W^{i\hspace{1.5mm}j}_{\phantom{i}\alpha \hspace{1.5mm} \beta} \partial _{i}
\partial _{j}\zeta ^{\beta } \partial_i \partial_j \zeta ,\zeta\> & = &
-\< W^{i\hspace{1.5mm}j}_{\phantom{i}\alpha \hspace{1.5mm} \beta} \partial _{i}
\partial _{j}\zeta ^{\beta } \partial_j u,\partial_i u\>+
\< \sigma(\nabla u) N, u\mid_{\partial \Omega}\>_{L^2(\partial \Omega)} 
\, , \vspace{1mm} \\
\< \nabla {\rm div}\,  \zeta, \zeta\> & = & 
-\< {\rm div}\, \zeta , {\rm div}\, \zeta \> + \< {\rm div}\,
 \zeta \, N,\zeta\>_{L^2(\partial \Omega)} \, .
\end{array}
$$
Just notice the novel
 boundary contribution arising in the last of these two terms, 
due to the fact that we now have $J(\zeta)=J(\BOne)=1$ over $\partial \Omega$, 
as opposed to ${\rm div} \, u\mid_{\partial \Omega} =0$ in the
analogous term in \S \ref{tde}, which then did not generate boundary 
contribution at all. The discretization of these two terms is straightforward.

The remaining term $\nabla q$ in (\ref{fif}) is dualized in a manner similar
to the procedure used in \S \ref{tde} for $\nabla h$, modulo one 
adjustment. Indeed, we first split $q$ as $q=q_0+ q(u)$, where $q(u)$ is the 
solution of 
$$
\begin{array}{rcl}
\Delta q(u) & = & \rho \,  {\rm trace}(\nabla u)^2 \, , \\
q(u)\mid_{\partial \Omega} & = & 0 \, .
\end{array}
$$
We have that $\nabla q = \nabla q_0 + \nabla q(u)$, and since 
$\nabla q(u)$ is a known function of $u$, we push it onto the right side 
of (\ref{fif}). We then proceed with $\nabla q_0$, and obtain that
$$ 
\< \nabla q_0 , \zeta \> =\;-{\displaystyle
\langle W^{i\hspace{1.5mm}j}_{\phantom{i}\alpha \hspace{1.5mm} j} N^\alpha
 , \nu \circ \zeta  \rangle } \< N, \zeta \> 
  - \< q_0, {\rm div} \zeta \>\, ,
$$
discretizing the summands on the right as we did their alteregos (arising 
from  $\nabla h$) in \S \ref{tde}.
 
We let
$$
R(u(t),F(t))= -\nabla q(u) -\rho v(t)+ \int_0^t F(s) ds
+\gamma^2 \rho (\BOne + \int_0^t u(s) ds)\, .  
$$
Using the numerical solution $U(t)$, and the discretized pressure force 
$F$ in the right side of (\ref{e24f}), we obtain a numerical representation
$R(U(t),F(t))$ of $R(u(t), F(t))$ that lies in 
$L^2_{1,b}(K', \mb{R}^3) \subset L^2_1(K';\mb{R}^3)$. 
Subject to the constraints (\ref{dfibc}), 
the vector of coordinate functions of (\ref{dfif}) is the solution of the
linear system of equations   
$$
\left(\gamma^2 \rho \BOne -A_{\circ, \partial}\right)
 \left( \begin{array}{c}
d_{e'_\circ} \\ d_{e'_{\partial}} \\
d_{f'_\circ} \\ d_{f'_{\partial}}
\end{array}
\right)= 
\left( \begin{array}{c}
\< R(U(t),F(t)), W_{e_\circ}\>  \\ \< R(U(t),F(t)), W_{e_{\partial}}\> \\
\< R(U(t),F(t)), W_{f_\circ}\>  \\ \< R(U(t),F(t)), W_{f_{\partial}} \>
\end{array}
 \right) \, ,
$$
where $A_{\circ, \partial}$ is the matrix 
$$
A_{\circ , \partial}
= \left(  \begin{array}{cccc}
\mb{O}_{e_\circ,e'_\circ} &  \mb{O}_{e_\circ,e'_{\partial}} &
\mb{O}_{e_\circ,f'_\circ} &  A_{e_\circ,f'_{\partial}} \\
\mb{O}_{e_{\partial},e'_\circ} &
A_{e_{\partial},e'_{\partial} } &
A_{e_{\partial},f'_\circ} & A_{e_{\partial},f'_{\partial} } \\
\mb{O}_{f_\circ,e'_\circ} &  A^t_{e_\partial,f'_{\circ}} &
K_{f_\circ,f'_\circ} &  K_{f_\circ,f'_{\partial}} \\
A^t_{e_\circ,f'_{\partial}} & A^t_{e_{\partial},f'_{\partial} } &
K_{f_{\partial},f'_\circ} & A_{f_{\partial},f'_{\partial} }
\end{array}
\right) \, ,
$$
the blocks $K_{f_\circ,f'_\circ}$, $K_{f_\circ,f'_{\partial}}$, and
$K_{f_\partial,f'_\circ}$ defined as they were in \S \ref{tde} 
(with $\lambda_L=1$), and the remaining nonzero blocks defined by  
$$
\begin{array}{rcl}
A_{e_\circ ,f'_{\partial}} & = & \frac{1}{2}
\left(\<{\rm div}\, W_{f'_{\partial}} N,W_{e_{\circ}}\> -(\<
  W^{i\hspace{1.5mm}j}_{\phantom{i}\alpha \hspace{1.5mm} j} N^\alpha
 , \nu \circ W_{e_\circ}\rangle \< N, W_{f'_{\partial}} \>
+\< W^{i\hspace{1.5mm}j}_{\phantom{i}\alpha \hspace{1.5mm} j} N^\alpha
 , \nu \circ W_{f'_{\partial}} \> \< N, W_{e_{\circ}} \>)\right)\, , \vspace{1mm} \\

A_{e_\partial ,e'_{\partial}} & = & -\frac{1}{2}(\< 
  W^{i\hspace{1.5mm}j}_{\phantom{i}\alpha \hspace{1.5mm} j} N^\alpha
 , \nu \circ W_{e_\partial}\rangle \< N, W_{e'_{\partial}} \>
+\< W^{i\hspace{1.5mm}j}_{\phantom{i}\alpha \hspace{1.5mm} j} N^\alpha
 , \nu \circ W_{e'_{\partial}} \> \< N, W_{e_{\partial}} \>)\, , \vspace{1mm} \\

A_{e_\partial ,f'_{\circ}} & = &
\frac{1}{2} \<\sigma(\nabla W_{f'_\circ})N, W_{e_\partial}\>
+ \frac{1}{2}\<{\rm div}\, W_{f'_{\circ}} N,W_{e_{\partial}}\>
- \vspace{1mm} \\
& & \frac{1}{2}(\<
  W^{i\hspace{1.5mm}j}_{\phantom{i}\alpha \hspace{1.5mm} j} N^\alpha
 , \nu \circ W_{e_\partial}\rangle \< N, W_{f'_{\circ}} \>
+\< W^{i\hspace{1.5mm}j}_{\phantom{i}\alpha \hspace{1.5mm} j} N^\alpha
 , \nu \circ W_{f'_{\circ}} \> \< N, W_{e_{\partial}} \>)
\, ,\vspace{1mm} \\

A_{e_\partial ,f'_{\partial}} & = &
\frac{1}{2}  \<\sigma(\nabla W_{f'_\partial})N, W_{e_\partial}\>
+ \frac{1}{2}\<{\rm div}\, W_{f'_{\partial}} N,W_{e_{\partial}}\>
- \vspace{1mm} \\
& & \frac{1}{2}(\<
  W^{i\hspace{1.5mm}j}_{\phantom{i}\alpha \hspace{1.5mm} j} N^\alpha
 , \nu \circ W_{e_\partial}\rangle \< N, W_{f'_{\partial}} \>
+\< W^{i\hspace{1.5mm}j}_{\phantom{i}\alpha \hspace{1.5mm} j} N^\alpha
 , \nu \circ W_{f'_{\partial}} \> \< N, W_{e_{\partial}} \>)\, ,\vspace{1mm} \\ 
\end{array}
$$
the pairings on the right here in the sense of $L^2(\partial \Omega)$, and
$$
\begin{array}{rcl}
 A_{f_{\partial},f'_{\partial} }  & = &  
 -\< \partial_i W^\alpha_{f_{\partial}},
A^{\alpha \beta}_{ij} \partial_j W^\beta_{f'_{\partial} }\>
- \< \partial_i W^i_{f_\partial }, \partial_j W^j_{f'_\partial}\>\\ & &
+\frac{1}{2}(\< l_{L_i(D)}\partial_i W^i_{f_\partial}, \partial_j
W^j_{f'_\partial }\>
+\< l_{L_i(D)}\partial_i W^i_{f'_\partial}, \partial_j W^j_{f_\partial}\>)
\vspace{1mm} + 
B_A(W_{f_{\partial}}, W_{f_{\partial}'})\,  ,
\end{array}
$$
where the new boundary term $B_A(W_{f_{\partial}}, W_{f_{\partial}'})$
is given by the sum of $L^2(\partial \Omega)$ pairings
$$
\begin{array}{rcl}
B_A(W_{f_{\partial}}, W_{f_{\partial}'}) & = &
\< \sigma(\nabla W_{f_{\partial}'})N,W_{f_{\partial}}\> +\frac{1}{2}
( \<{\rm div}\, W_{f'_{\partial}} N,W_{f_{\partial}}\>
+\<{\rm div}\, W_{f_{\partial}} N,W_{f'_{\partial}}\>) \vspace{1mm}\\
& & 
-
\frac{1}{2}(\<
  W^{i\hspace{1.5mm}j}_{\phantom{i}\alpha \hspace{1.5mm} j} N^\alpha
 , \nu \circ W_{f_\partial}\rangle \< N, W_{f'_{\partial}} \>
+\< W^{i\hspace{1.5mm}j}_{\phantom{i}\alpha \hspace{1.5mm} j} N^\alpha
 , \nu \circ W_{f'_{\partial}} \> \< N, W_{f_{\partial}} \>) \, .
\end{array}
$$

By construction, $Z(0,x)=x$, and $\dot{Z}(0,x)=\mc{W}(x)$, 
$\mc{W}$ the discretization 
in $L^2_1(K';\mb{R}^3)$ of the divergence-free initial velocity $w$. 
In general, $Z(t,\Omega) \neq \Omega$ and 
$Z(t,\partial \Omega) \neq \partial \Omega$ 
when $t>0$.

\section{Simulation results} 
\label{s6}
We take five of the negative eigenvalues of 
$-(\rho I_{\circ, \partial})^{-1} K_{\circ, \partial}$, 
$-(\rho I_{\circ, \partial^B})^{-1} K_{\circ, \partial^B}$ for the
matrices in the systems (\ref{e24}) and (\ref{e24f}), respectively, as
indicated below, and describe the ensuing coarse and fine resonance 
vibration patterns of $\Omega$ by depicting the nodal points of these 
waves over the portion of the boundary opposite to that where the 
external sinusoidal force $F=F_f$ hits it. 

\begin{figure}[H]
\includegraphics[height=3.0in]{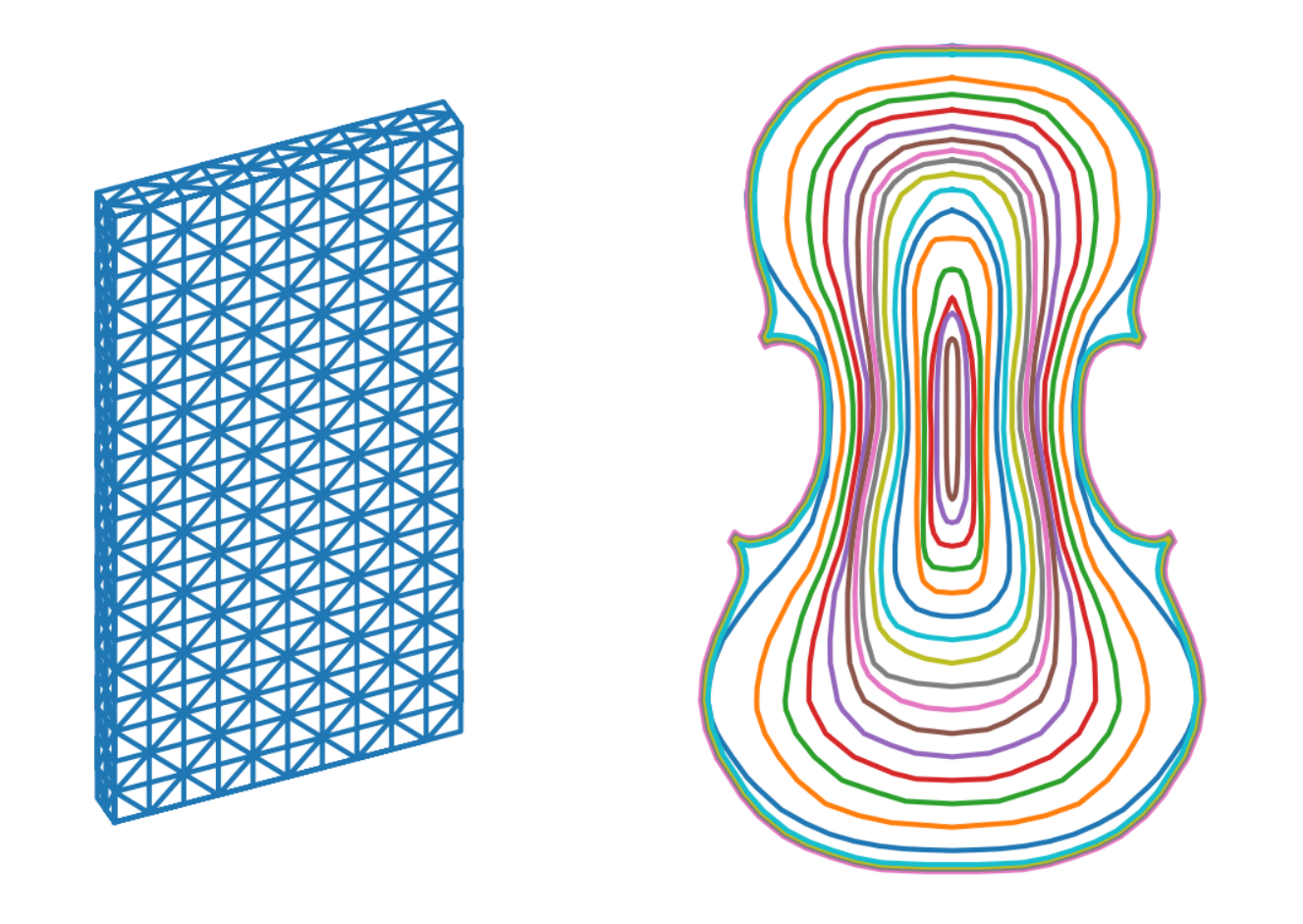}
\caption{Left: Slab with its triangulation, as used here. Right: Viotti
top plate view (without the $f$ holes) and sixteen level sets.}
\label{NF1}
\end{figure}

In our experiments, we repeat the geometries considered in \cite{sim}:
\begin{enumerate}
\item In the first of the experiments, $\Omega$ is a slab of  
$\text{$10$ cm}\times \text{$1$ cm} \times \text{$20$ cm}$, or 
thinner versions of depth $0.5$ cm and $0.25$ cm, 
respectively, Fig. \ref{NF1} left.
\item For the second of our experiments, $\Omega$ has the geometry of the top 
plate of the classic Viotti violin, as per \cite{tstrad}, without its
$f$-holes, Fig. \ref{NF1} right.
\end{enumerate}

The slabs of different depths serve to check the effects that 
flatness of the boundary of $\Omega$, and rescaling in the thin direction, 
have on the vibration patterns. 

\subsection{Computational complexity}
The slab of $\text{$10$ cm}\times \text{$1$ cm} \times \text{$20$ cm}$ 
is subdivided into $400$ regular 
blocks of size $\text{$1$ cm}\times \text{$0.5$ cm} \times \text{$1$ cm}$ each, 
and triangulated accordingly, with the blocks given the standard subdivision 
into five tetrahedrons each. The triangulation $K$ so obtained has 
$693$ vertices, $3,212$ edges, 
$4,520$ faces and $2,000$ tetrahedrons; $522$ of the vertices, $1,560$ of the
edges, and $1,040$ of the faces, are on the boundary. The first barycentric 
subdivision $K'$ contains $10,425$ vertices, $61,544$ edges 
$99,120$ faces, and $48,000$ tetrahedrons; $3,122$ of the vertices,
 $9,360$ of the edges, and $6,240$ of the faces, are on the boundary.
The two other thinner versions of the slab have triangulations with the same 
number of elements, merely scaling the depth accordingly. 
The aspect ratio of their tetrahedrons are a half, and a quarter of the aspect
ratios that they have in the first of the triangulations, respectively. 
In this case, (\ref{e24}) is a system of 
$160,664$ equations in $160,664$ unknowns, an increase of  
$62\%$ in size of the corresponding system treated in \cite{sim}.

For $\Omega$ the top plate of the Viotti violin without its $f$-holes,
the triangulation is more complex. This  
$\Omega$ can be inscribed into a rectangular box of
$35.4\,  {\rm cm} \times 20.8 \,  {\rm cm}$; it is curved with
thickness that varies nonuniformly, 
ranging from a lowest of $0.21$ cm to a largest of 
$0.36$ cm.    
In order to obtain a reasonable resolution for this type of thickness and
curvature, away from the edge of the plate, we subdivide $\Omega$ into blocks of
size $\text{$0.5$ cm}\times \text{$0.5d$ cm} \times \text{$0.5$ cm}$,
where $d$ is an average thickness of the plate at points where the block is
located. We proceed similarly at the edge, but the first dimension of the 
blocks we consider there is taken to be nonuniform, out of necessity.
 All together, it takes 4,608 of these blocks to fill $\Omega$, and
the triangulation $K$ that we then derive contains 7,296 vertices, 
35,283 edges, 51,028 faces, and 23,040 tetrahedrons;
$5,118$ of the vertices, $15,348$ of the edges,
and 10,232 of the faces, are on the boundary.
The first barycentric subdivision of this triangulation has 116,647 vertices,
699,294 edges, 1,135,608 faces, and 552,960 tetrahedrons; 
$30,698$ of the vertices, $92,088$ of the edges,
  and 61,392 of the faces, are located on the boundary.
The system (\ref{e24}) has now 
$1,834,902$ equations in $1,834,902$ unknowns, again, an increase of  
$62\%$ in size of the corresponding system in \cite{sim} (the agreement
in the increase between this body and the slab is due to the 
comparable aspect ratios of their geometries and triangulating elements).

All of these $\Omega$s are topologically $3$d balls with a $2$d sphere boundary.
The number of elements in their triangulations $K$ and $K'$ satisfy 
the combinatorial Euler characteristic identity of the $3$d ball,   
and the number of boundary elements satisfy the combinatorial Euler 
characteristic identity of the $2$d sphere.

In either case, as we traverse the body across the two blocks separating 
the bounding surface in the thin direction, the barycentric triangulations
used contain 444 faces, with 35 intermediate faces separating pair of 
boundary faces on opposite ends, which when included yield a total 37 faces 
altogether traversed in going from one side to the other. These two 
blocks contain $382$ edges, so there are just about $32$ edges in between the 
two said faces, and so in going from one side to the other in the thin
direction, we cross nearly $69$ triangulation elements that are involved in the
expansion of the discretized approximate solution $U$.  
This provides an adequate resolution for the numerical approximation to 
accurately capture the true nature of the vibrating wave in these regions, 
and to propagate in all transversal directions it goes, or feel the propagating
waves that are passing by. Refinements of the 
triangulations would increase the number of elements in going from side to side,
and improve the accuracy of the numerical solution $U$, but such 
would lead to a complexity that is out of the scope of the author's current 
computational resources.

As for the additional local boundary conditions that are satisfied by the fine 
vibration waves approximations, of all the $2 \times 6$ 
subblocks of (\ref{e27}) that the faces in ${\mc F}_{\partial }$ generate, 
of the 1,040 of them for the slab, exactly eleven 
have rank one, and of the 10,232 of them for the Viotti plate, exactly
one has this property \cite{sim}. As pointed out 
in \cite{sim}, the differences between the numbers of null rows in the
system (\ref{e27}) for the slab and 
Viotti plate reflect the $\mb{Z}/2 \times \mb{Z}/2 \times \mb{Z}/2$ symmetry 
of the triangulation of the former, as opposed to the nonsymmetric nature of 
the triangulation of the latter, whose thickness varies in a 
nonuniform manner throughout the body. Notice that
the exterior normal to a boundary face, and its normal as an 
oriented two simplex, are not necessarily the same.

We find the matrices $I$ and $K$ of the systems (\ref{e24}), (\ref{e24f}),  
by executing python code written for this purpose. The code is structured into
four major modules. The first two are common to all geometries and topologies, 
and the last two are body specific, dealing with the coarse and fine waves, 
respectively.
(We have improved the object-oriented conception used in \cite{sim},
the organization of the code making it straightforward the incorporation 
of bodies with other physical and topological properties by a mere insertion
of the new appropriate modules in the right places.)
The processing of the results, and their graphical 
display, is carried out by some additional python code written for 
the purpose, the graphical component of it built on top of the PyLab 
standard library module.

We executed the code on a 2.4GHz Intel Core i7 processor, with 8GB 1600 MHz 
of memory. The generation of the new blocks in $I$ and $K$ associated to 
edges only require a CPU time that is of the same order as that taken 
to generate the blocks corresponding to faces only, reported in \cite{sim}. 
For the number of elements involved, as well as the sparsity of the blocks 
for edges, are of the same order of magnitude as those for the faces.
The hardest new part arises when generating the cross block 
corresponding to interior edges and interior faces,  
given the number of elements involved, and the lower sparsity that this 
off-diagonal subblock has in relation to the others (for the Viotti plate, 
the sparsity score of the $I_{\circ,\partial}$ matrix 
of the system (\ref{e24}) here was $0.999990$; in contrast, the 
sparsity score of the corresponding matrix in \cite{sim} was $0.999994$). 
Overall, the hardware used handles well the complexities of the problem, 
but the size of it for the Viotti plate already demands the generation of the 
matrices in stages, parallelizing the calculations on the basis of the 
linear ordering of the simplices in $K'$, and finding them employing a 
still reasonable amount of CPU time. A true technological problem 
arises in the calculation of the eigenvectors of (\ref{e24}) and
(\ref{e24f}) for this plate, which with the \verb+eigsh+ ARPACK routine we 
use for the purpose, requires the use of a bigger serial system than any 
we have had available to us (see \S \ref{s31} below).

The effective cross sectional area of the Viotti plate is about $2.88$ times  
the effective cross sectional area of the slab, and this produces a 
difference of one order of magnitude in the number of elements, edges and 
faces, that we must use in treating one or the other. Should we need to
treat a thin $\Omega$ with cross sectional area $2.88$ times larger than that
of the Viotti plate, the processor(s) needed for the purpose
should be capable of handling matrices of size $10^7\times 10^7$, which
is well within the reach of today's computers.
Our method scales, and thus, seems suitable for the treatment of several 
imaginable problems of interest in the current technological 
environment.

\subsection{Elastic constants}
The components of the tensor of elastic constants $W$ in (\ref{ten2}) are
those given in (\ref{ten3}). For this tensor $W$, the constants 
$l_1, l_2, l_3$ in the weighted divergence (\ref{wedi}) are
$$
\begin{array}{rcl}
l_1 & = & 1.932758876\cdot 10^9\, ,  \\
l_2 & = & 1.488135884\cdot 10^9\, ,  \\
l_3 & = & 5.884378014 \cdot 10^9 \, . 
\end{array}
$$
As indicated earlier, we take $\rho=360 {\rm kg/m^3}$, and $\lambda_L=1$.
These are the values of the constants used in our simulations. 
The density value implies that the masses of the slabs are
$18$g, $9$g, and $4.5$g in decreasing order of their thickness, 
respectively, and that the mass of the Viotti plate is approximately 
$59.1$g, very close to the actual mass of many violin plates of this size 
currently in existence (the Messiah, for instance). Just for the 
the perspective of the reader, we observe that the density 
of aluminum, $2,700$ ${\rm kg/m^3}$, is $7.5$ times $\rho$.


\subsection{Simulations} 
We analyze first the divergences of the eigenvector wave modes 
of the systems (\ref{e24}) and (\ref{e24f}), respectively, and then 
do the analysis of the resonance waves that they produce. 
 
\subsubsection{The divergence of the coarse and fine normalized 
eigenvector solutions} \label{s31}
The eigenvalues and eigenvectors of the homogeneous systems
associated to (\ref{e24}) and (\ref{e24f}) are generated using the 
ARPACK routine \verb+eigsh+ in shift-invert mode, with their corresponding 
matrix parameters $\rho I$, and $K$, respectively.  This computes the 
solutions $(\lambda, c)$ of the system 
$$
\rho I c = \lambda K c \, .
$$
With \verb+sigma=-1/(2 pi f)^2+, and \verb+which='LM'+ passed onto 
\verb+eigsh+, we execute the routine for a frequency 
\verb+f+ any of 
$80$, $147$, $222$, $304$, and $349$ Hz, respectively. In each case, 
the routine returns pairs \verb+-((2 \pi f_r)^2, c^{f_r})+ of eigenvalue and 
eigenvector of $(\rho I)^{-1}K$, where \verb+f_r+ is the eigenvalue
of the matrix closest to the inputted \verb+f+, in magnitude. We then
consider the corresponding {\it normalized} eigenvector wave solution for each, 
the coarse, and fine systems, respectively. 

For the Viotti plate simulations, the execution of $\verb+eigsh+$ above to find
the eigenvectors fails with the error 
\verb+Can't expand MemType 1: jcol 1801648+, barely missing the end.
It is for this reason that we adjust these simulations to use only 
the face elements in the triangulation, and the corresponding altered 
equation. These are the results we present below for the Viotti plate.
Notice that the altered equation is not the same as that considered in
\cite{sim}, because the definition of the matrix $K$ here differs from what it
was then.

In the case of (\ref{e24}), we denote by $f_{coarse}^r$ and 
$c_s^{f^r_{coarse}}$ the pair \verb+(f_r,c^{f_r})+ produced for the inputted
\verb+f+, and let  
$$
U_{f^f_{coarse}}=e^{2\pi f^r_{coarse} t i}\sum_{W_s\in 
\mc{B}_{L^2_1(K';\mb{R}^3)}} c^{f^r_{coarse}}_s W_s
$$
be the normalized coarse eigenvector solution in $L^2_1(K';\mb{R}^3)$
that results. In the case of the system (\ref{e24f}), we 
proceed likewise, and let 
$$
U_{f^f_{fine}}=e^{2\pi f^r_{fine} t i}\sum_{W_s \in
\mc{B}_{L^2_{1,b}(K';\mb{R}^3)}} 
c^{f^r_{fine}}_s W_s
$$
be the normalized fine eigenvector solution in $L^2_{1,b}(K';\mb{R}^3)$ that 
results from the returned pair \verb+(f_r,c^{f_r})+ for the given \verb+f+.

We study the divergence of any of these normalized eigenvector solutions  
$U$ by computing the flux 
$$
\int_{\partial \Omega} U \cdot n \, d\sigma 
$$
through the boundary of the body at time $t=0$. The results are displayed in
Table 3.
\medskip

\begin{table}
\begin{tabular}{|c|r|r|r|r|r|}
\hline \hline
Body & $f$\phantom{3} & $f^r_{coarse}$\phantom{333} & flux $U_{f^f_{coarse}}$ &
 $f^r_{fine}$\phantom{333} & flux $U_{f^f_{fine}}$ \\ \hline \hline
\multirow{5}{*}{Slab 1.0} & 80 & 79.98300620 & -0.0000249730 
&  80.01751279 & 0.0000009708 \\
& 147 & 146.90402861 & -0.0001222083 
& 146.01597845 & 0.0000663559 \\
& 222 & 221.93558743 &  0.0003161676 
& 220.45892181 &-0.0057717711 \\
& 304 & 304.03517536 & -0.0023043236 
& 304.01774121 &-0.0000023969 \\
& 349 & 348.94594451 & -0.0000171669 
& 348.96922189 & 0.0000051104 \\ \hline
\multirow{5}{*}{Slab 0.5} & 80 & 79.44465641 & -0.0041503791 
&  79.44465959 & -0.0009114188 \\
 & 147 & 145.97955476 &-0.0083605924 
& 145.97955400 & -0.0020807749 \\ 
 & 222 & 220.45898203 &-0.0026676193   
& 220.46140952 & -0.0002427300 \\ 
 & 304 & 301.88969619 &-0.0076778693 
& 301.88974794 & -0.0016486990 \\ 
 & 349 & 349.09691350 &-0.0000202374 
& 349.05496747 & -0.0000142132 \\ \hline  
\multirow{5}{*}{Slab 0.25} & 80 & 79.44465624 &-0.0031694079 
& 79.44465547 & -0.0063644096 \\
 & 147 & 145.97955476 &-0.0059265938 
&145.97955448 & -0.0063411167 \\
 & 222 & 220.45898620 &-0.0018536358 
&221.99574371 &  0.0000703494 \\
 & 304 & 301.88969891 &-0.0050515305 
&301.88969505 & -0.0041960055 \\
 & 349 & 348.90773121 &-0.0000748571 
&349.09787762 &  0.0000857175 \\ \hline
\multirow{5}{*}{Viotti plate} & 80 & 79.99818471 &-0.0087776487 
       &  79.99972869 & -0.0070498420  \\
 & 147 & 146.99513939 & -0.0168672080 
       & 146.99453429 &  0.0173401990  \\
 & 222 & 222.00259999 &  0.0112724697 
       & 221.99478910 & -0.0005874142 \\
 & 304 & 303.99865189 & -0.0053818663 
       & 304.00167789 & -0.0094033063 \\
 & 349 & 348.99889409 & -0.0036936708 
       & 349.00660655 &  0.0029127309 \\
\hline \hline
\end{tabular}
\smallskip

\caption{Initial fluxes for the coarse and fine eigenvector 
solutions. The results for the Viotti plate are derived using the face elements
only.}
\end{table}
\medskip

We repeat this experiment just for the slab of depth $0.5$cm using   
the value $\lambda_L =3.101757591\cdot 10^9$, the 
average of $l_1$, $l_2$, and $l_3$ above, instead of $1$. Then, the spatial 
terms in equation (\ref{vp2}), which are all of order two, have coefficients
that are of the same order of magnitude, and for certain frequencies, the 
modifying weak term 
$-\lambda_L \< {\rm div} \, u, {\rm div}\, u\>$ could drag the other two 
microlocally towards an operator that is not elliptic, producing a non 
physical wave result that could be detected by some numerical inconsistency,
for instance, a negative ``norm'' for the eigenmode that is to be     
normalized in order to compute its flux.
This does not happen here at any of the five frequencies used in our 
simulations, but could in principle occur for others that remain unexplored.  
The complete results are displayed in Table 4.
\medskip
 
\begin{table}
\begin{tabular}{|c|r|r|r|r|r|}
\hline \hline
& $f$\phantom{3} & $f^r_{coarse}$\phantom{333} & 
flux $U_{f^f_{coarse}}$ & $f^r_{fine}$\phantom{333} & 
flux $U_{f^f_{fine}}$ \\ \hline \hline
\multirow{5}{*}{ \mbox{$
\begin{array}{c}
\text{Slab 0.5} \\ \lambda_L=3.1 \cdot 10^9 
\end{array}
$}
} & 80 &  79.44465548 & 0.0085159297   
 & 79.44465548 & 0.0172430719 \\
 & 147 & 146.91260194 & 0.0001067346  
 & 145.97955453 & 0.0188712978 \\
 & 222 & 222.22211502 &-0.0000724956  
 & 220.45891928 &-0.0196273836 \\
 & 304 & 301.88969256 &-0.0080961614  
 & 301.88969160 & 0.0256947693 \\
 & 349 & 346.57731074 &-0.0082506915  
 & 346.57731061 & -0.0316121099 \\ 
\hline \hline
\end{tabular}
\smallskip

\caption{Slab 05: Initial fluxes for the coarse and fine 
eigenvector solutions of the equation with $\lambda_L=3.101757591\cdot 10^9$.}
\end{table}

\subsubsection{Coarse and fine resonance waves} 
We subject the body to an external sinusoidal pressure wave of the form 
$$
\vec{F}^{ext}=\vec{F}_0 \sin({\bf k}\cdot {\bf x} \mp \omega t) 
=\vec{F}_0 \sin({\bf k}\cdot {\bf x})\cos({\omega t})\mp
\vec{F}_0 \cos({\bf k}\cdot {\bf x})\sin({\omega t}) 
$$
that travels in the appropriate direction for it to hit the bottom 
of the slab, or belly of the Viotti plate, first. If $k=|{\bf k}|$ is the
magnitude of the wave vector, we use $v=\omega/k=343$ m/sec,  
the speed of sound in dry air at $20\mbox{}^\circ$C. The source of the
wave is placed at a distance of $62$cm from the body, along a line 
that passes through the height-width plane of the body perpendicularly at 
the half way point of both, its height and width. This external force 
$\vec{F}^{ext}$ induces a force on the body, which we denote by
$\vec{F}_{\omega}$.

The coarse (\ref{e24}) and fine (\ref{e24f}) systems are considered with the  
nonhomogeneous force term $F=\vec{F}_{\omega}$,
and with trivial initial data. 
(These are the coarse and fine discrete versions of the initial value problem 
for (\ref{vp}) when this system is viewed as the second order equation 
(\ref{vp2}).)  The nonhomogeneous terms in 
these systems are vectors of the form
$$
\vec{C}_1 \cos{(\omega t)}\mp \vec{C}_2 \sin{(\omega t)}\, ,
$$
where $\vec{C}_1$ and $\vec{C}_2$ are time independent vector fields on the
body. 
If $\omega =2\pi f$ for $f$ any of the frequency values 
$80$, $147$, $222$, $304$, and $349$ Hz, respectively, 
we let $\omega_{f_r} =
2\pi f_r$ be the closest eigenvalue to $\omega$ that the matrix 
$(\rho I)^{-1}K$  
has, $\rho I$ and $K$ the matrices of the system in consideration.
The resonance wave that this vibration mode produces is given by  
\begin{equation} \label{resw}
W_{f_r}=\frac{1}{\omega^2_{f_r}-\omega^2}\sum_{W_s}
\left(c^{W_s}_1\left(-\cos{(\omega_{f_r}t)}+
\cos{(\omega t)}\right)\pm c^{W_s}_2 \left( \frac{\omega}{\omega_{f_r}}
\sin{(\omega_{f_r}t)}-\sin{(\omega t)}\right)\right)W_s \, ,
\end{equation}
where, for each $j=1,2$, the vector $\vec{c}_j =(c_j^{W_s})$ is the solution 
to the linear system of equations  
$$
-((2\pi f)^2 \rho I +K)\vec{c}_j = \vec{C}_j \, . 
$$
The summation is over the basis elements $W_s$ of $L^2_1(K';\mb{R}^3)$ and
$L^2_{1,b}(K';\mb{R}^3)$ for the coarse and fine systems, respectively.

We generate the vector $\vec{C}_j$ as a function of $\vec{F}_\omega$, and 
the geometry of the body, and then solve the system of equations above for 
$\vec{c}_j$ 
using the \verb+scipy.sparse.linalg+ routine \verb+spsolve+, with the 
appropriate parameters. 

We compute the values of the resonance wave at the barycenter, and vertices, 
of the boundary faces on the side of the body opposite to the 
incoming external wave $\vec{F}^{ext}$. 
There are $3,661$ such points for the slabs, and 
$42,069$ for the Viotti plate. We do these calculations at the equally spaced 
times $t_j=j\frac{2\pi}{10\omega}$, $j=1 , \ldots, 10$, corresponding to a 
full cycle of the external wave. For each $t_j$, we find the maximum
${\rm max}_{t_j}=\max{\{\| W_{f_r}(t_j)\|}\}$ and minimum 
${\rm min}_{t_j}=\min{\{\| W_{f_r}(t_j)}\|\}$ of the set of norms of the 
resonance wave at the indicated points, and with 
$\delta_{t_j}=\frac{1}{10}({\rm max}_{t_j}-{\rm min}_{t_j})$,  
any of the said points is considered to be nodal at time $t_j$ if the norm 
of the resonance wave solution at the point is no larger than 
${\rm min}_{t_j}+ c_\Omega \delta_{t_j}$, with $c_\Omega=0.05$ for the slabs,
and $c_\Omega=0.00002$ for the Viotti plate, respectively. 
We then declare a point to be nodal if it is nodal at all the $t_j$s. 
These are the points that we display. 

The results for the coarse and fine resonance waves are depicted
in Figs. \ref{f2}-\ref{f5} below. In each case, we indicate
the value of $c_\Omega$ that is being used to define a point as nodal.    

The triangulations with simplices of best aspect ratios are those for
the slab of depth $1$cm, and the Viotti plate. We  restrict our attention
to them, and study the change in the resulting resonance pattern produced
by taking into consideration now the six modes with eigenvalues
$\omega_{f_{r_j}}=2\pi f_{r_j}$,
$j=1,\ldots, 6$, closest to $\omega=2\pi f$, as opposed to the single
closest one, as above. The resonance wave is then a sum of six terms as
in (\ref{resw}), one per each of the $f_{r_j}$s. We show the results of
these simulation for the slab plate resonating at the frequency
$f=147$Hz only, Fig. \ref{f6}.
The five cases depicted correspond to nodal points defined as above for
values of $c_\Omega=0.1, 0.05, 0.03, 0.02$, and $0.005$,
respectively. We bypass showing the corresponding image for the Viotti plate
because of the difficulties encountered with \verb+eigsh+ to find
its resonance vibration modes using edges and faces together. The image
we obtain using the face elements only compares very well to its analogue
\cite[Figure 7]{sim}, the equation solved here being better when studying 
the resonance patterns over a wider range of vibrations, not just the small 
ones back then.

Finally, we look at the divergence of the resonance waves
above. We denote by $W_{f_{coarse}^f}$ and $W_{f_{fine}^f}$ the normalized
coarse and fine resonance wave solutions associated to the pair $(f,f_r)$ for
(\ref{e24}), and (\ref{e24f}), respectively. As these waves start with trivial
initial conditions, we quantify the extent to which our algorithm maintains
the divergence free condition on them throughout time by evaluating their
fluxes over the boundary at the time $t_j$, in each case, where
the quotient ${\rm max}_{t_j}/{\rm min}_{t_j}$ is the largest.
The normalization performed on the resonance waves
makes the results independent of the magnitude of $\vec{F}_0$
in the external wave $\vec{F}^{ext}$ that induces the resonance. The said
$t_j$ (at which value the fluxes are being evaluated) turns out to be
a nontrivial function of $\Omega$ and $f$, as opposed to the analogous
situation in \cite{sim}, where it was equal to $t_7$ always.
The results are listed in Table 5.

\begin{figure}[H]
\includegraphics[height=6.00in,angle=-90]{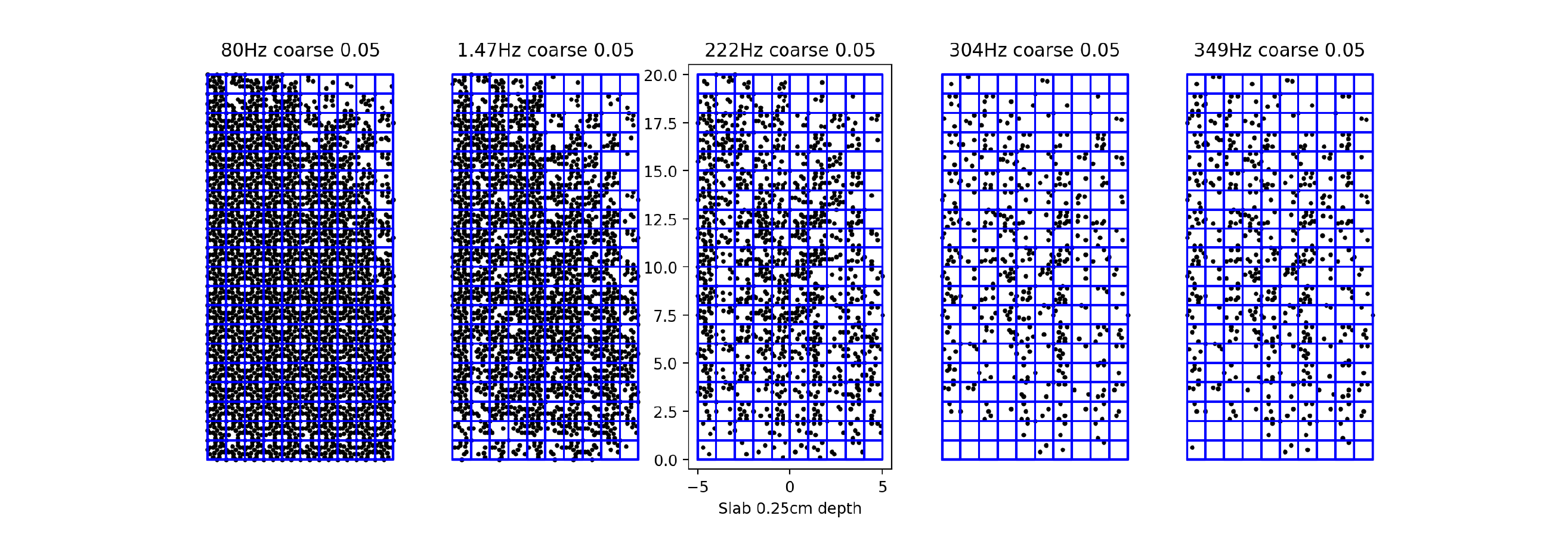}

\includegraphics[height=6.00in,angle=-90]{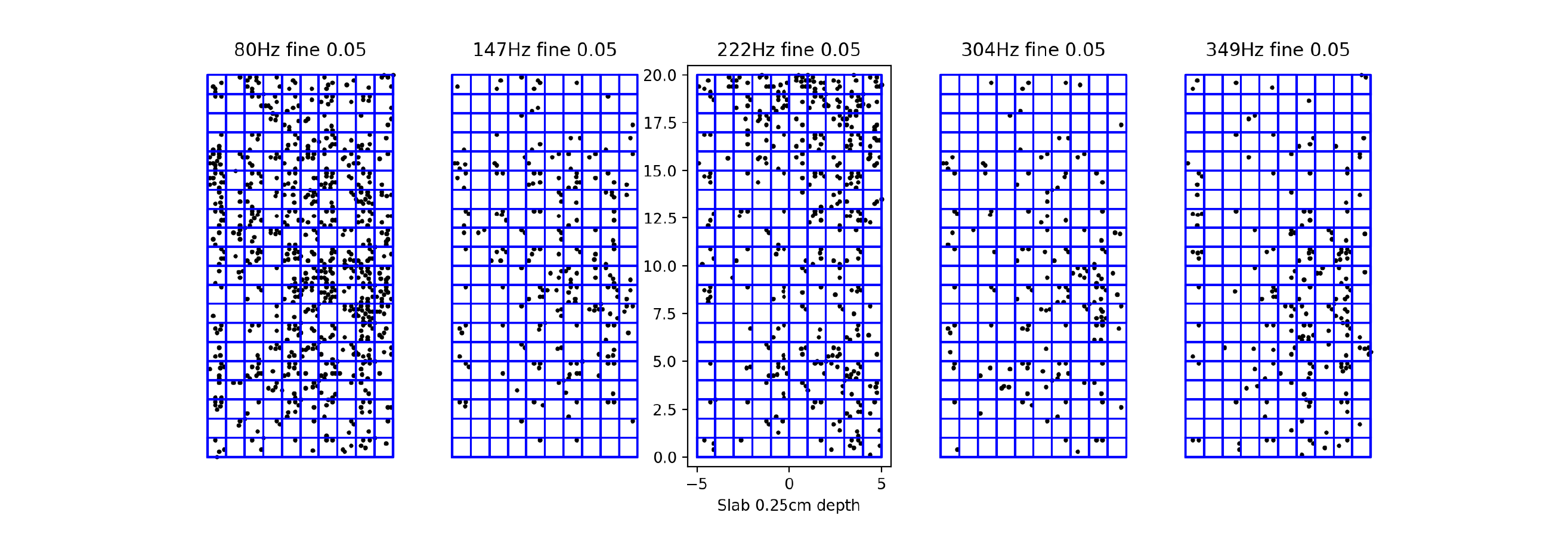}
\smallskip

\caption{Slab $10\times 0.25 \times 20$: Nodal points of the resonance waves
at $f$Hz arising from the slab mode of vibration of frequency $f_r$ closest
to $f$. The notion of nodal point is defined by $c_\Omega=0.05$.}
\label{f2}
\end{figure}
\medskip

\subsection{Discussion of the results: Validation}
The elastic constant values in Tables 1 and 2
for the Engelmann Spruce had been collected from wood with 
a 12\% moisture content \cite{gwk}, which is not very
conducive to the creation of good vibration patterns, and although 
the density value we employ in our simulations is close to the density of 
actual wood used in the making of violin plates, if compared to the
original, our model for the Viotti plate is dull, and vibrates poorly.
(D. Caron, a renowned american luthier, makes his violin plates using
wood with a moisture content in the range of 2\%-5\% \cite{ca}; this content
is lowered when the plate is coated with varnish, which adds
mass and absorbs some moisture as it dries. The addition of the
varnish changes the flexibility by stiffening the plate across the grain,
in effect, creating an almost functionally graded plate, though not 
isotropic.) In spite of this drawback, the vibration patterns of the
bodies in our simulations are very close to what they would be in practice
under those nonideal circumstances.

\begin{figure}[H]
\includegraphics[height=6.00in,angle=-90]{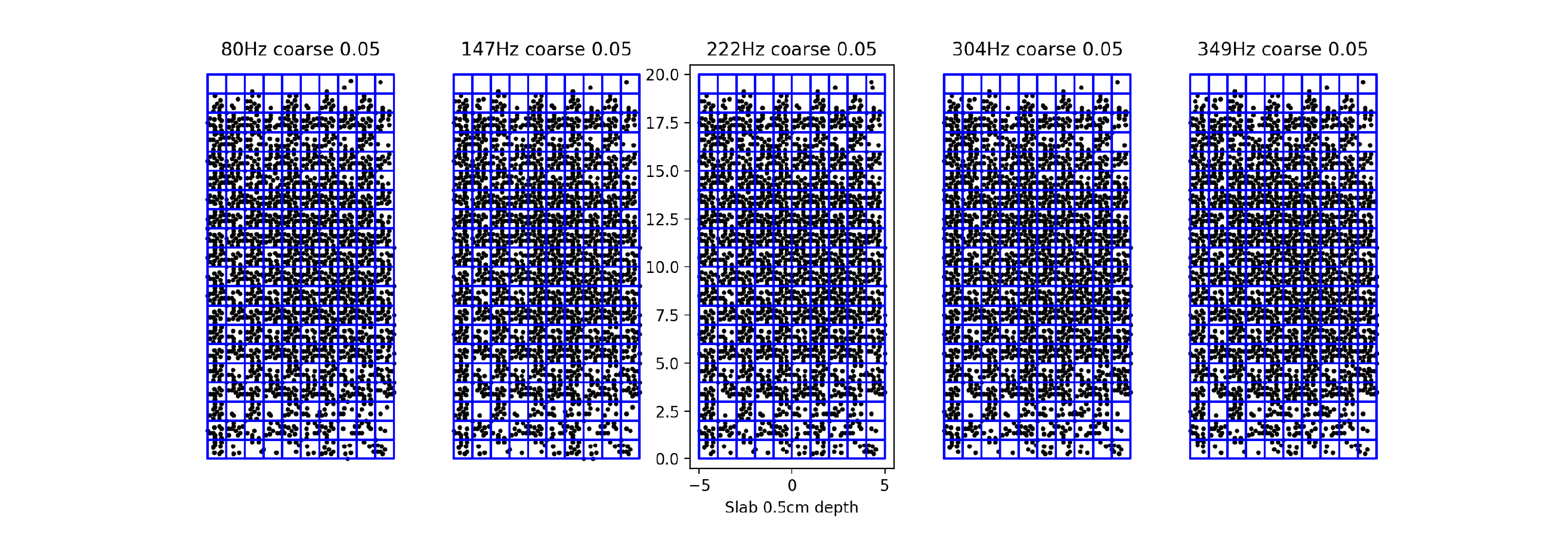}

\includegraphics[height=6.00in,angle=-90]{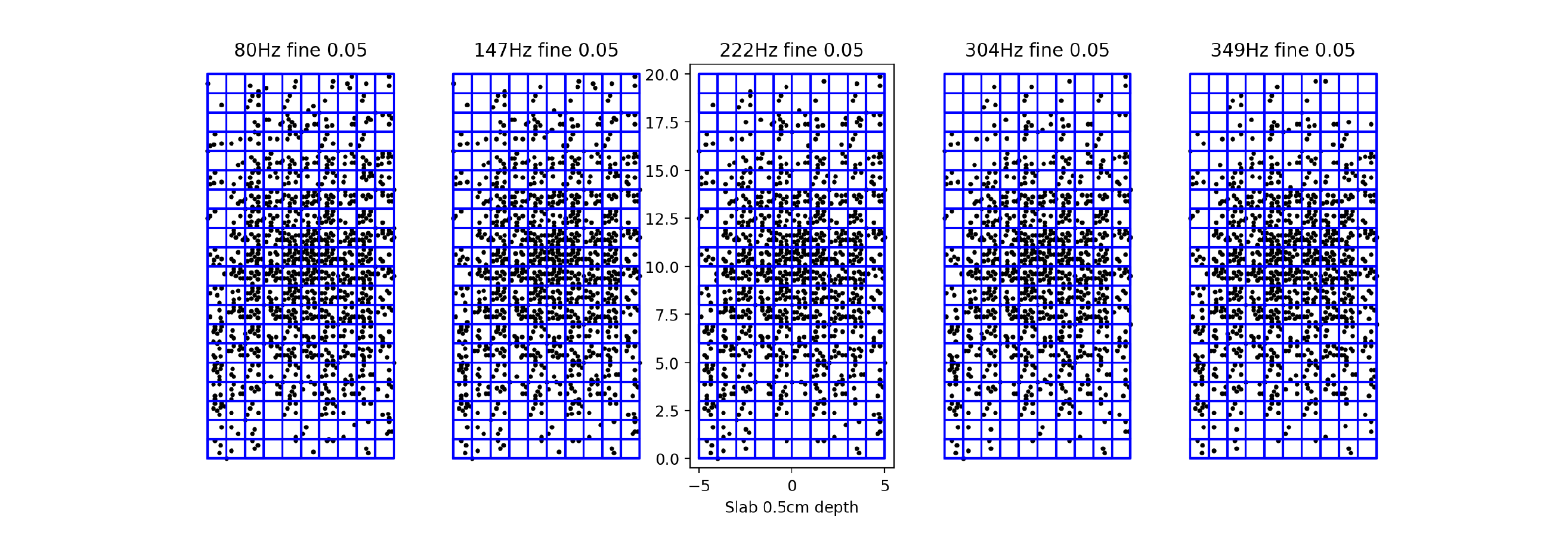}
\smallskip

\caption{Slab $10\times 0.5 \times 20$: Nodal points of the resonance waves
at $f$Hz arising
from the slab mode of vibration of frequency $f_r$ closest to $f$. The notion
of nodal point is defined by $c_\Omega=0.05$.}
\label{f3}
\end{figure}

\begin{figure}[H]
\includegraphics[height=6.00in,angle=-90]{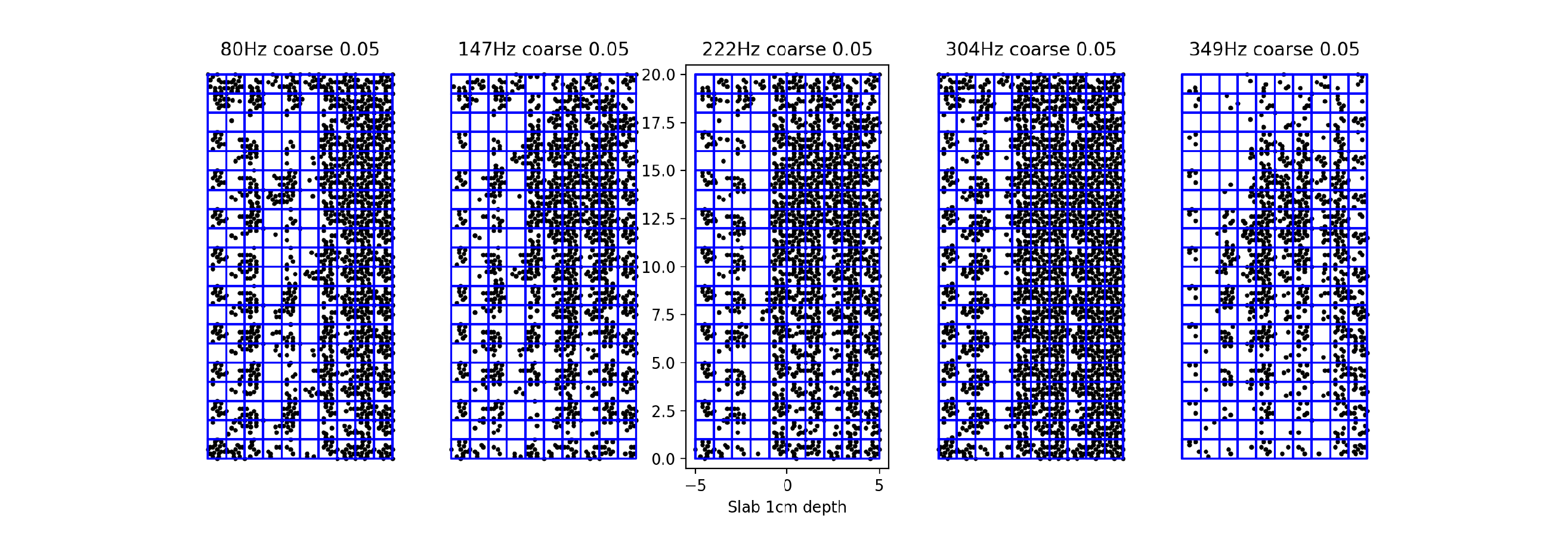}

\includegraphics[height=6.00in,angle=-90]{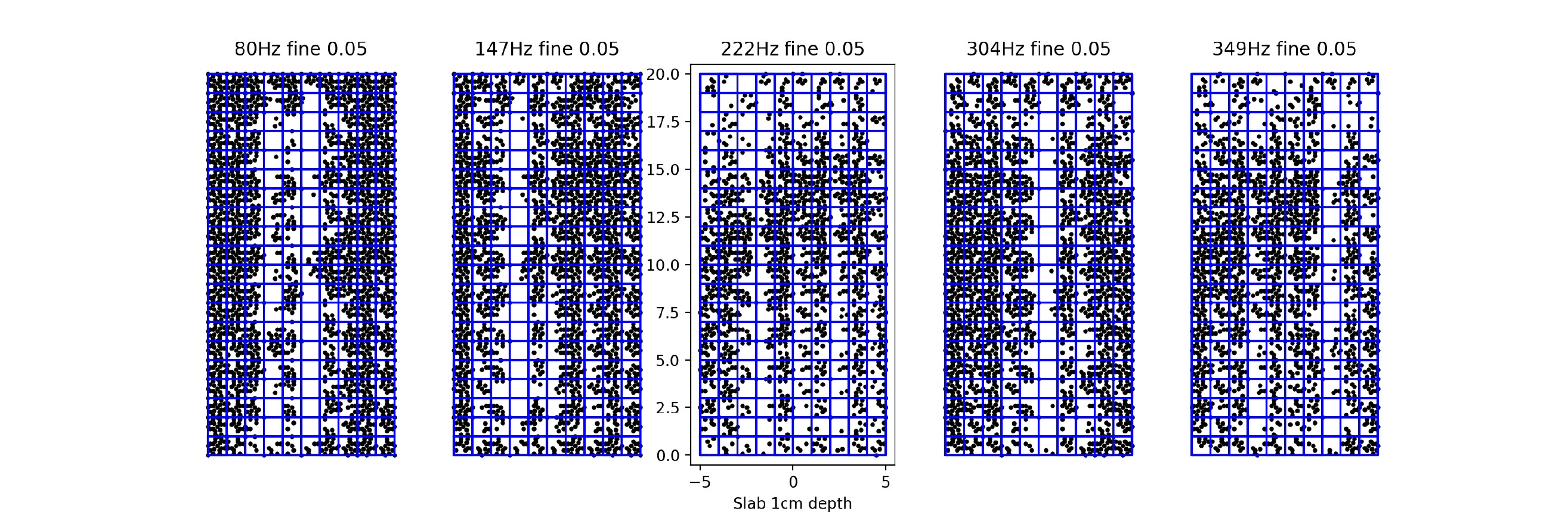}
\smallskip

\caption{Slab $10\times 1.0 \times 20$: Nodal points of the resonance waves
at $f$Hz arising
from the slab mode of vibration of frequency $f_r$ closest to $f$. The notion
of nodal point is defined by $c_\Omega=0.05$.}
\label{f4}
\end{figure}

\begin{figure}[H]
\includegraphics[height=6.00in,angle=-90]{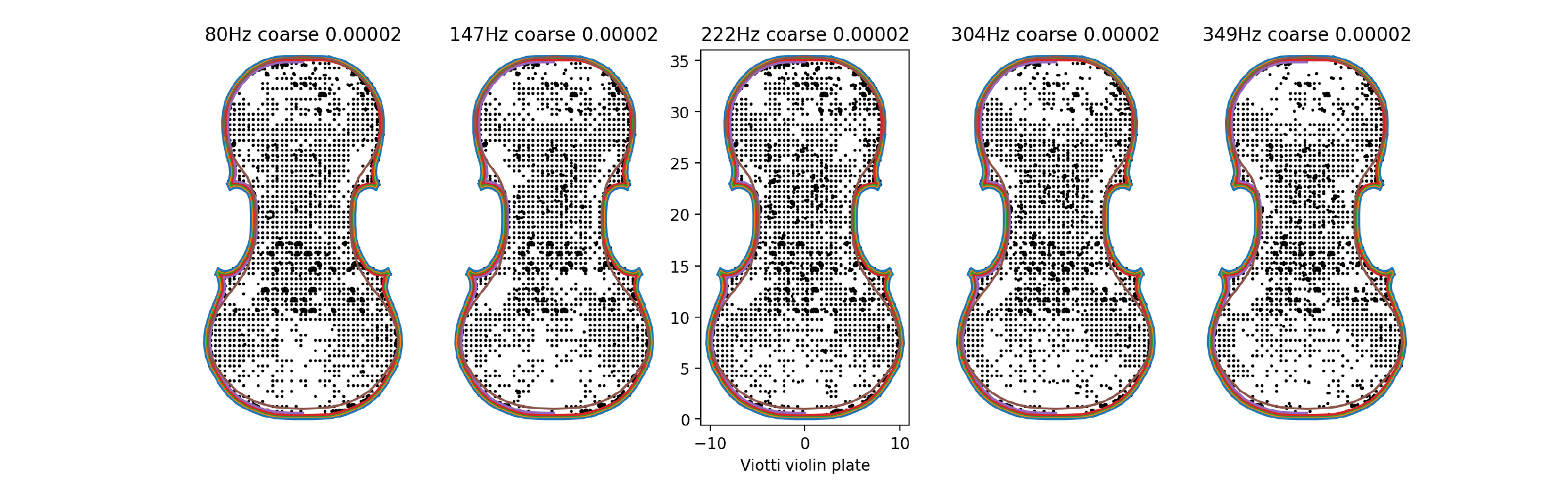}

\includegraphics[height=6.00in,angle=-90]{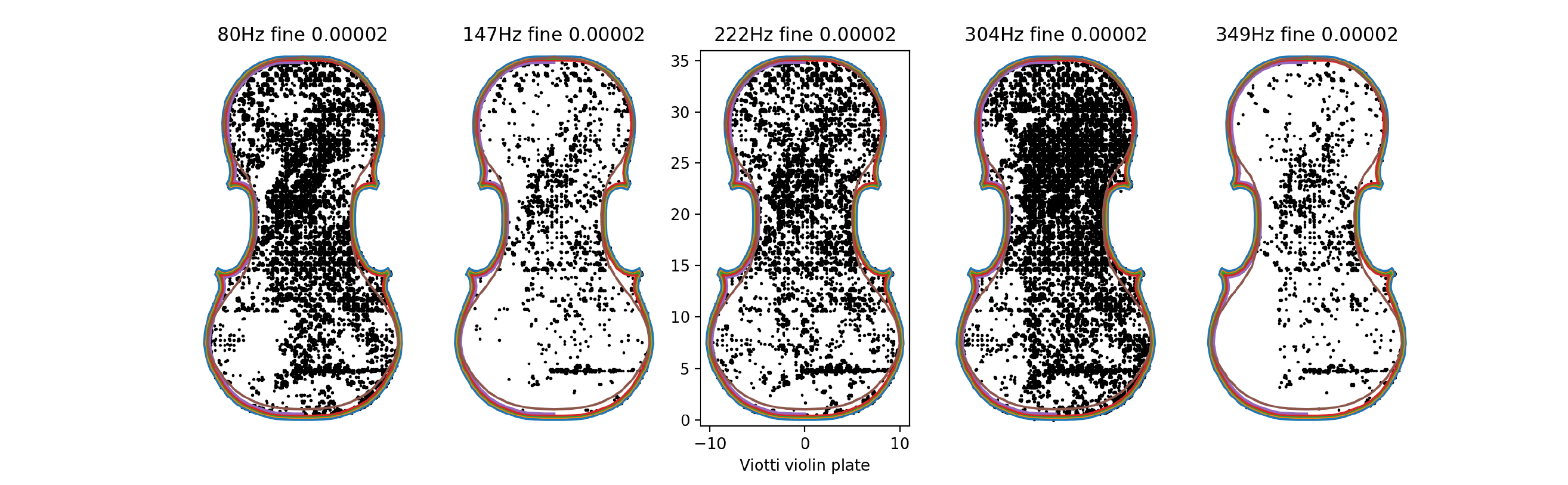}
\smallskip
\caption{Viotti plate: Nodal points of the resonance waves at $f$Hz arising
from the plate mode of vibration of frequency $f_r$ closest to $f$ (results
derived using the face elements only). The
notion of nodal point is defined by $c_\Omega=0.00002$.}
\label{f5}
\end{figure}

\begin{figure}[H]
\includegraphics[height=6.00in,angle=-90]{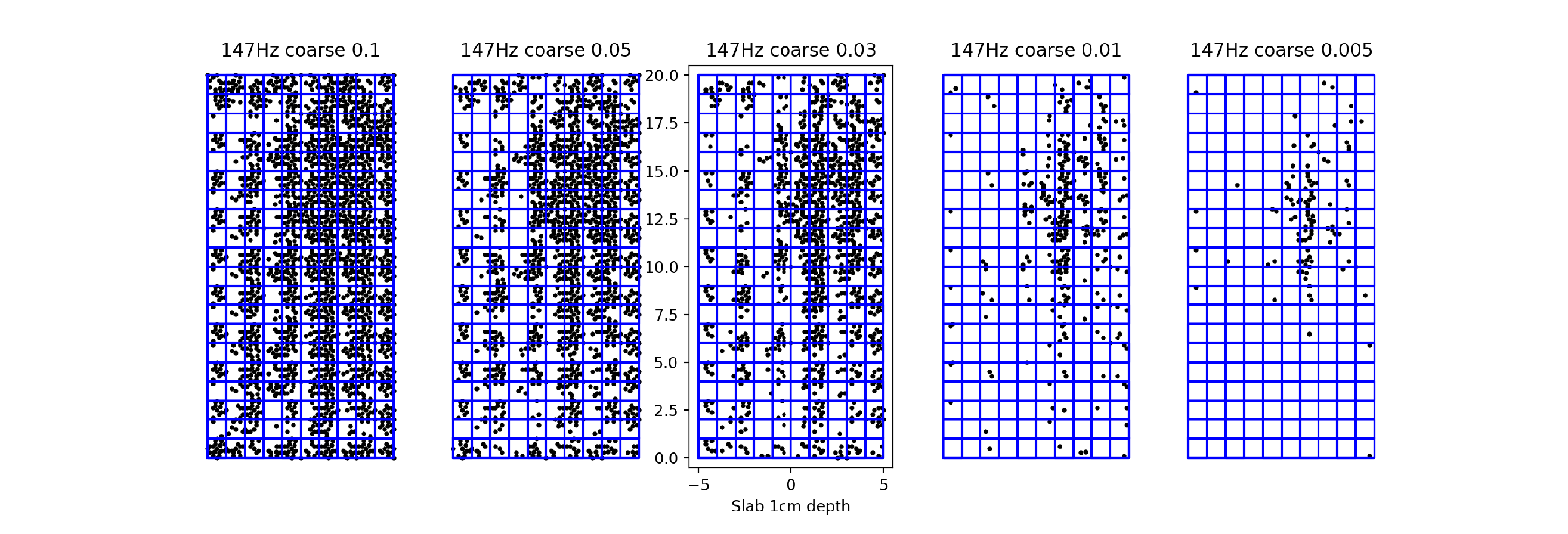}

\includegraphics[height=6.00in,angle=-90]{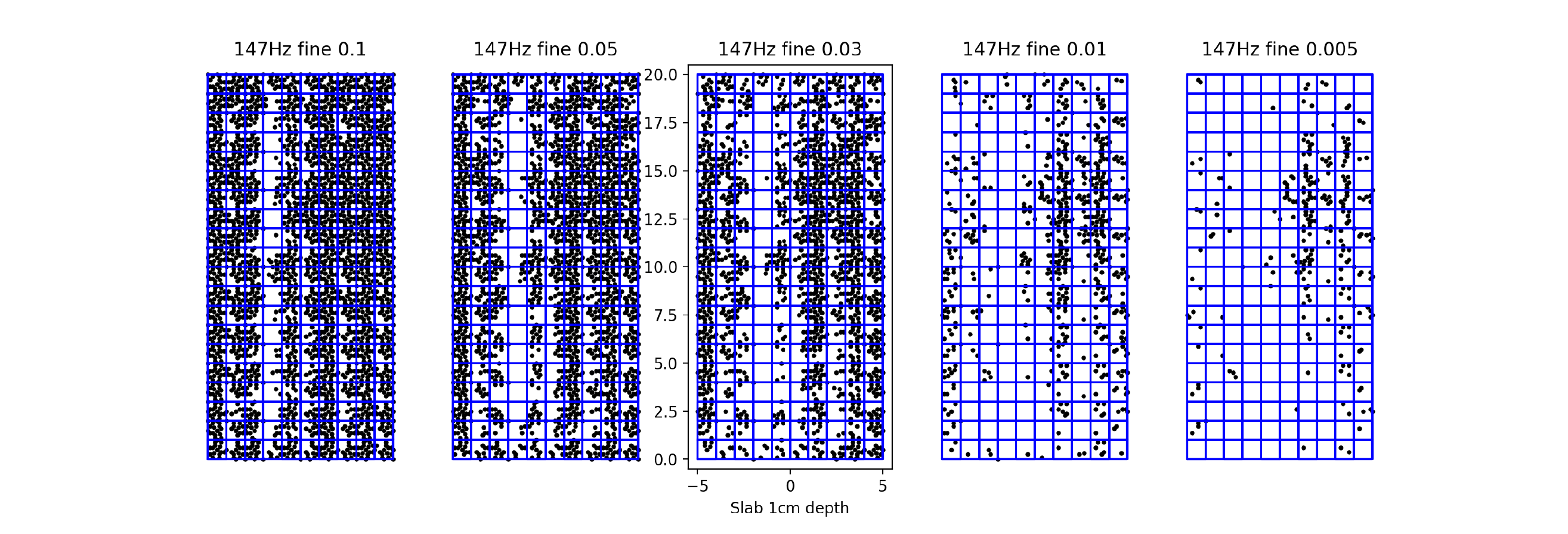}
\smallskip

\caption{Nodal points of the resonance wave at $f=147$Hz arising from
the six closest mode of vibrations of the plate to this frequency $f$.
From left to right, the notion of nodal point is
defined by $c_\Omega=0.1$, $0.05$, $0.03$, $0.01$,  and $0.005$,
respectively.}
\label{f6}
\end{figure}

Figs. 2-4, and 6 depict resonance patters in remarkable agreement with 
holographic images of standing sound waves, propagating in a body under 
experimental conditions comparable to those in our simulations 
(see, for instance, \cite[Figs. 4, 6, 7]{ruph}), showing that they are  
what they are supposed to be. The 
corresponding resonance patterns in Figs. \ref{f4} and \ref{f6}  
exhibit just a better definition of the details if we take into consideration
the six vibrations modes closest to $147$Hz, as opposed to the single 
closest one, but not an actual change in the pattern, and in Fig. \ref{f6},  
as the notion of a nodal point becomes stricter by decreasing the value of
$c_\Omega$, clearer details in the resulting patterns emerge. 
The simulation at $147$Hz for the Viotti plate in Fig. \ref{f5}, which is 
derived using the face elements only, points quite closely towards the 
image by holographic interferometry of the mode 2 of a top violin plate 
in \cite[Fig. on p. 177]{hu}, and improves on those shown in
\cite[Figs. 6, 7]{sim}. In spite of the differing conditions between our 
simulations and the cited holographic experiments, we
take the favorable comparison as a validation of
our results. We elaborate on a few additional details of our simulations.

\medskip

For the two most geometrically degenerate bodies in our simulations, 
at a given frequency, the number of nodal points of the coarse waves is 
larger than the number for the fine waves. As the frequencies increase, the 
ratios of these numbers are 2.08, 2.15, 2.37, 2.51, and 2.41 for the slab 
of depth $0.5$cm, and 10.62, 10.88, 3.23, 3.13 and 2.88 for the slab of 
depth $0.25$cm, respectively. 
These patterns breaks altogether for the slab of depth 
1cm, and for the Viotti plate, the two bodies in our simulations with
triangulations of simplices of best aspect ratios. For comparison
purposes, the said quotients for the thicker slab are
0.73, 0.71, 1.13, 1.14, and 0.60, while for the Viotti plate they are
0.33, 1.32, 0.49, 0.22, and 1.48, respectively.

A close look at the number of nodal points for the slabs in our simulations 
indicate a limitation to the effect of selective thinning of the plate, a 
practice among luthiers that is commonly believed to lower or increase 
frequency modes if carried out in regions of high or low curvature, 
respectively.  For the thickest of the slabs the number of fine nodal points
varies somewhat sinusoidally with $f$, oscillating about 1,990. For the
slab of depth $05$cm, the numbers go as 1,076, 1,050, 965, 928, and
945, , while for the thinnest slab they go as 276, 215, 381, 158, and 188,
respectively.
At some point, the effect of additional thinning acts differently
on the various vibration modes, counterbalancing the strength they have
relative to each other.
\medskip

The curvature of the violin plate is an important factor in the relative
magnitude of the amplitude of the resonance wave near $f=349$Hz compared
to the amplitude of the resonance wave near $f=80$Hz. The
quotient of these amplitudes for the Viotti plate is of the order
of $10$, while for the slab of depth $1$cm is of the order of
$2\cdot 10^3$. (The brightness of this vibration mode is reported to
be of importance in prejudging the good quality of an assembled instrument
with the plate in it \cite{hu}.) Of course, the comparison is weak in that
the simulation for the Viotti plate was carried out using only the face 
elements, however many. But the resonance pattern that we have so obtained,
Fig. \ref{f5}, in spite of the poor quality of the wood used, is 
remarkably good. 

\section{Concluding remarks: Computational physics merging with
computational topology}
\label{s7}
If $\eta$ is a solution of (\ref{eq1}), (\ref{eq2}), (\ref{eq3}), we have that
 ${\rm div}_{\eta}\dot{\eta}=0$, and the linearized equations of motion
about $(\eta,\dot{\eta})$ have as solution a curve $(u(t),\dot{u}(t))$ with
$u$ a divergence-free vector field that, while well-defined, yields a curve
in the Abelian group $H^2(\Omega;\mb{R})$. Solving for $u$ then is
fraught with the problems inherent to maintaining the closed
divergence-free property that the initial condition satisfies, hard
theoretically, and harder computationally.

The linearized equations about a point $(\eta ,\dot{\eta })$ other
than $(\eta(t)x=x,0)$ produces an operator that is hyperbolic on the tangent
space at $\eta$ of the submanifold defined by (\ref{eq1}). The dependence
of the space on $\eta$ makes it difficult to prove the well-posedness of its
associated Cauchy problem in the usual manner. We enlarge
the space to prevent this problem, but this in turn forces us to
modify the equation in order to maintain the ellipticity of the spatial
part of the linearized operator, and have its hyperbolicity on the
larger domain. This is what we accomplish when we modify (\ref{eqs}),
and look at (\ref{new2}) instead,
the latter an equation whose volume preserving solutions
solve (\ref{eqs}) as well. Technical issues aside,
we then show that the fixed point $(\eta,\dot{\eta})$ of the
contraction mapping principle that we devise solves the modified equation,
and is volume preserving, so this fixed point solves (\ref{eqs}) as well,
which is what we wanted.

\begin{table}
\begin{tabular}{|c|r|r|r|r|r|}
\hline \hline
Body & $f$\phantom{3} & $f^r_{coarse}$\phantom{333} & flux $W_{f^f_{coarse}}$ & $f^r_{fine}$\phantom{333} & flux $W_{f^f_{fine}}$ \\ \hline \hline
\multirow{5}{*}{Slab 1.0} & 80 & 79.98300620 & -0.0045986774 
&  80.01751279 & -0.0000092823 \\
& 147 & 146.90402861 & -0.0022899377 
& 146.01597845 & -0.0000271994 \\
& 222 & 221.93558743 & -0.0074495371 
& 220.45892181 & -0.0000039583 \\
& 304 & 304.03517536 &  0.0006411380 
& 304.01774121 &  0.0000089606 \\
& 349 & 348.94594451 & -0.0141437966 
& 348.96922189 &  0.0000717342 \\ \hline
\multirow{5}{*}{Slab 0.5} & 80 & 79.44465641 &  0.0000233879 
&  79.44465959 & 0.0000114983 \\
& 147 & 145.97955476 &  0.0000740003 
& 145.97955400 & 0.0000103866 \\ 
& 222 & 220.45898203 &  0.0001275103 
& 220.46140952 & 0.0000111290 \\ 
& 304 & 301.88969619 &  0.0003156497 
& 301.88974794 & 0.0000175093 \\ 
& 349 & 349.09691350 &  0.0004423489 
& 349.05496747 & 0.0000111168 \\ \hline
\multirow{5}{*}{Slab 0.25} & 80 & 79.44465624 & -0.0000157439 
&  79.44465547 &   0.0000021530 \\
& 147 & 145.97955476 &  0.0000016388 
& 145.97955448 &   0.0000103168 \\
& 222 & 220.45898620 & -0.0000246448 
& 221.99574371 &   0.0000710929 \\
& 304 & 301.88969891 &  0.0002215706 
& 301.88969505 &   0.0000297836 \\
& 349 & 348.90773121 &  0.0003536785 
& 349.09787762 &  -0.0000881610 \\ \hline
\multirow{5}{*}{Viotti plate} & 80 
      &  79.99818471 & -0.0189303539 &  79.99972869 &  0.0035747785 \\
& 147 & 146.99513939 & -0.0575648428 & 146.99453429 & -0.0977355895 \\
& 222 & 222.00259999 & -0.0189054794 & 221.99478910 & -0.0118013729 \\ 
& 304 & 303.99865189 & -0.0072224699 & 304.00167789 &  0.0000439681 \\
& 349 & 348.99889409 & -0.0036936708 & 349.00660655 & -0.0922473868 \\ 
\hline \hline
\end{tabular}
\smallskip

\caption{Fluxes of the normalized coarse and fine resonance 
waves computed each at the value of $t_j$ maximizing 
${\rm max}_{t_j}/{\rm min}_{t_j}$ (results for Viotti plate 
derived using the face elements only).}
\end{table}

Having the existence settled, we can then turn our attention to the finding 
of a good numerical approximation to it.
As it is impossible to maintain a closed condition when you solve 
numerically any equation, we give up on attempting to compute 
the volume preserving diffeomorphism solution of the equation, and use 
the iterates of the Newton scheme in the proof of its existence as 
the candidates to approximate it instead. 
The velocity fields of the diffeomorphisms in the iteration are not 
divergence-free per se, but close to one. Such a condition, a matter 
of satisfying an equation that locally involves finitely many of the 
coefficients in the discretized unknown, is left to be regulated by the 
equation itself, and if holding initially, it should be satisfied at later 
times within a small margin of error also due to the well-posedness of the 
problem solved, thus ensuring that the nonlocal incompressible condition 
of the diffeomorphism solution stays within a striking level of tolerance. 
Thus, we can accurately captured it numerically if we use discretizing spaces 
for them that are natural relative to the complexities of the problem, 
and that encode into them the algebraic topology of the actual solution 
velocity curve in the cohomology group above that this velocity represents,
a global property. This is where we enter a triangulation $K$ of $\Omega$, and 
all of the Whitney forms of the simplices of its barycentric subdivision $K'$.
Through them we resolve both of these last issues in earnest.

For any $k$, $0\leq k \leq n$, the Whitney forms of the simplices in $K'$ 
of degree $k$ are a finite basis of the $k$th cochain group of $K'$. Since the 
simplices are all contractible, we may define locally a notion of 
$*$ operator, and use the $*$ of the Whitney forms of degree $n-k$ 
simplices to generate a dual cochain group of $k$-simplices in $K'$.
Preserving its algebraic topology properties, 
we discretize any $k$-form in the direct sum of these two groups, 
a set-up that is closely related to the dual block decomposition of $K'$,
and the proof of the Poincar\'e duality $H^{n-k}(K';\mb{R})\cong H_{k}(K';
\mb{R})$. Using the individual groups for this purpose would lead to the loss 
of algebraic topology information, although some is still preserved.     
With that resolved, the issue of a proper resolution for an accurate 
numerical approximation of the said $k$-form is addressed by controlling the
aspect ratios of all the simplices in $K$, making of this a relatively
uniform number throughout the entire polytope $|K|$. A minimum
lower bound on this number ensures a resolution good enough to have 
the form accurately described by its numerical approximation. This is what
we have done in our problem, where $n=3$, and $k=1$, except for the additional
fact that instead of looking at the $k$ form, we have worked equivalently
with its metric dual vector field.

The large systems of equations that need to be solved, however sparse, are 
computationally challenging, and this is 
a difficulty that we must face even after resolving that of determining the 
systems themselves. The judicious choice of discretizing spaces minimizes 
these difficulties. Since we encode in them the initial conditions of the
waves we seek, with  the algebraic topological properties that they have,
we can then allow the well-posedness of the equations solved for to play its 
role in keeping the numerical solutions derived physically accurate in time.
  
The use of Whitney forms associated to edges for numerical purposes of sort
is well established already \cite{rao, asv, boss}, but often they have been
employed to discretize physical quantities that truly represent cohomology 
classes in degree two, rather than one, and for which the use of the 
Whitney forms associated to faces would be more natural instead 
\cite{sim}, the direct sum of the forms for edges and faces described here 
a better candidate to work with, as a discretizing space, in general.
Historically, as computers started to be developed, functions were 
approximated by their values at finitely many points chosen 
uniformly over their domains, the number of these points limited by
the memory capacity, rather 
than expanding them as linear combinations of the elements in the partition 
of unity given by the Whitney forms associated to vertices, or tetrahedrons, 
individually, or together as a direct sum. That to a great extent 
this has remained so is kind of surprising, given the common use of 
orthonormal basis of $L^2$ functions derived from a starting wavelet 
\cite{dau1,dau2,mall}, and which has appeared frequently in the numerical 
analysis of signals in a noisy background, the orthonormality of the basis 
ensuring that in the representation of the 
sequence of images as a discrete set of observation processes, the noise is 
no more correlated than it was in the original images \cite{hasi}. The
space ${\rm Char}(K')$ serves perfectly well in that capacity also, an
orthogonal partition of unity basis of $L^2$ functions generated by the Haar 
wavelet supported on a model top dimensional simplex. But the sum
${\rm Fun}(K') \oplus {\rm Char}(K')$ would be better for general purposes,
as we had seen in \S \ref{s4}, the natural space to discretize functions as
a sum of their components in the image of the Laplace operator, and 
their orthogonal complements. 

Regardless, the use of these forms is being increasingly advocated 
today in numerical analysis \cite{boss2}, and have regained some of 
the momentum they had at the time of Whitney (who had used them for several 
``computations'' already), and shortly thereafter \cite{do}. 
If their properties are exploited well, 
they serve to capture the essential algebraic structure underlying the 
problems under consideration, and, with a minimal number of computational 
elements, we may use that fact to derive accurate approximations of solutions
to problems with large intrinsic complexities. 
The geometric content of a triangulation is a powerful tool to use 
to compute polytope quantities of physical significance \cite{gsi}, 
a tool through which we may control the local geometry of the problem. 
But the power of this tool truly comes to fruition if we include in the 
considerations 
its algebraic content as well, since then we can use it to tie up, and control, 
both the local, and global geometric and topological properties of it.
 The serial computer power remains a 
hindrance only to the size of the problem that you can treat, but that is just
a technological issue. Parallelizing the calculations will take you as far as 
that may go, but if and when you need to proceed to do so, you know 
beforehand that the problem you are solving for is well posed. 

As we indicated earlier, we can also study the 
vibrational patterns of elastic bodies with corners, under mild 
assumptions on their stored energy functions, or even bodies where part of the 
boundary is fixed while the rest is free to move. These are all variations
of the theme treated here. What is of more importance is the fact that
the computation of any time evolving tensor of physical 
significance can be treated in this manner also. We just need to
add a minor part to the set-up above concerning the Whitney forms. For
a simplex $\sigma_k=[x_{p_0},x_{p_1}, \ldots, x_{p_k}]$, the Whitney form 
associated to it is the complete alternation of the tensor
$\tau_{\sigma_k}=x_{p_0} dx_{p_1} \otimes \cdots \otimes dx_{p_k}$. In general,
in order to treat evolving tensors that are not purely contravariant in a
manner similar to the treatment of the tensor in our problem here, we would 
have to use the complete symmetrization of all the $\tau_{\sigma_k}$s as well,
in parallel with the use of all of the Whitney forms, as above. The collection
of all of these local elements associated to the simplices of $K'$ 
is the optimal space where to track, at the discrete level, a weak 
version of the tensor.

\end{document}